\documentclass[12pt]{article}
\usepackage{amsmath, amsthm, amssymb,amscd, verbatim,hyperref,mathrsfs}
\theoremstyle{plain}
\newtheorem{thm}{Theorem}[section]
\newtheorem{thmn}{Theorem}
\newtheorem{lem}[thm]{Lemma}
\newtheorem{cor}[thm]{Corollary}

\newtheorem{prop}[thm]{Proposition}

\theoremstyle{definition}

\newtheorem{nota}[thmn]{Notation}
\theoremstyle{remark}
\newtheorem{rem}[thm]{Remark}
\newtheorem*{rema}{{\bf Remark 6.A}}
\newtheorem*{remb}{{\bf Remark 4.1.A}}
\newtheorem*{remc}{{\bf Remark 4.2.B}}

\newcommand{\nc}{\newcommand} 
\nc{\hb}{\mathbb} 
\nc{\M}{\mathcal} 
\nc{\mf}{\mathfrak}
\nc{\mbf}{\mathbf}
\nc{\DMO}{\DeclareMathOperator}

\newbox\noforkbox \newdimen\forklinewidth
\setbox0\hbox{$\textstyle\smile$}
\setbox1\hbox to \wd0{\hfil\vrule width \forklinewidth depth-2pt
 height 10pt \hfil}
\setbox\noforkbox\hbox{\lower 2pt\box1\lower 2pt\box0\relax}
\def\anchor{\mathop{\copy\noforkbox}\limits}
 
\setbox0\hbox{$\textstyle\smile$}
\setbox1\hbox to \wd0{\hfil{\sl /\/}\hfil}
\setbox2\hbox to \wd0{\hfil\vrule height 10pt depth -2pt width
               \forklinewidth\hfil}
\newbox\doesforkbox
\setbox\doesforkbox\hbox{\box1 \lower 2pt\box2\lower2pt\box0\relax}
\def\nanchor{\mathop{\copy\doesforkbox}\limits}
 
\nc{\cA}{{\M A}} \nc{\cB}{{\M B}} \nc{\cC}{{\M C}} \nc{\cD}{{\M D}}
\nc{\cE}{{\M E}} \nc{\cF}{{\M F}} \nc{\cG}{{\M G}} \nc{\cH}{{\M H}}
\nc{\cI}{{\M I}} \nc{\cJ}{{\M J}} \nc{\cK}{{\M K}} \nc{\cL}{{\M L}}
\nc{\cM}{{\M M}} \nc{\cN}{{\M N}} \nc{\cO}{{\M O}} \nc{\cP}{{\M P}}
\nc{\cQ}{{\M Q}} \nc{\cR}{{\M R}} \nc{\cS}{{\M S}} \nc{\cT}{{\M T}}
\nc{\cU}{{\M U}} \nc{\cV}{{\M V}} \nc{\cW}{{\M W}} \nc{\cX}{{\M X}}
\nc{\cY}{{\M Y}} \nc{\cZ}{{\M Z}}
\nc{\Aa}{{\hb A}} \nc{\Cc}{{\hb C}} \nc{\Gg}{{\hb G}}
\nc{\Nn}{{\hb N}} \nc{\Pp}{{\hb P}} 
\nc{\Qq}{{\hb Q}} \nc{\Rr}{{\hb R}} \nc{\Zz}{{\hb Z}}
\nc{\mfa}{{\mf a}} \nc{\mfb}{{\mf b}} \nc{\mfk}{{\mf k}}
\nc{\mfm}{{\mf m}} \nc{\mfp}{{\mf p}} \nc{\mfq}{{\mf q}}
\nc{\mfr}{{\mf r}}
\nc{\fP}{{\mf P}}
\DMO*{\trdeg}{td}
\DMO*{\spec}{Spec}
\DMO*{\fork}{\nanchor}
\DMO*{\dnf}{\anchor}
\DMO{\RU}{RU}
\DMO{\deter}{det}
\DMO{\RM}{RM}
\DMO{\RC}{RC}
\DMO{\Real}{Re}
\DMO{\Imag}{Im}
\DMO{\tr}{tr}
\DMO{\qc}{QC}
\DMO{\Hu}{Hull}
\DMO{\leg}{length}
\DMO{\area}{area}
\DMO{\dia}{diameter}
\DMO{\iso}{Iso}
\DMO{\dis}{dist}
\DMO{\grad}{grad}
\DMO{\vol}{volume}
\DMO{\gra}{grad}
\DMO{\hd}{nbhd}
\DMO{\dv}{div}
\DMO{\Psl}{PSL}
\nc{\Mb}{\mathfrak^{2b/\delta}_{K_x}}
\nc{\Ma}{\mathfrak^{2a/\delta}_{K_x}}
\nc{\dif}{\mathrm{d}}
\nc{\G}{\Gamma}
\nc{\g}{\gamma}
\nc{\D}{\nabla}
\nc{\p}{\partial}
\nc{\DD}{\Delta^2}
\nc{\pp}{\partial^2} 
\nc{\de}{\delta}
\nc{\td}[2]{\trdeg{({#1}/{#2})}}
\nc{\dtd}[2]{\trdeg_{\delta}{({#1}/{#2})}}
\nc{\dspec}[1]{\spec_{\delta}{#1}}
\nc{\ddim}[1]{\dimen_{\delta}{#1}}
\nc{\gens}[1]{\langle {#1} \rangle}        
\nc{\gen}[2]{ {#1} \langle {#2} \rangle } 
\nc{\form}{\Omega}
\nc{\set}[1]{\left\{ {#1} \right\}}
\nc{\mr}{\hat}
\nc{\pr}{\partial}
\nc{\bc}[3]{\cB^{#1}({#2},{#3})=B^{#1}_{#2}(C^{#2}_{#3}(t)+B^{#1}_{#3}(C^{#2}_{#3}(t))} 
\nc{\tuple}[2]{{#1},\ldots,{#2}} \nc{\ptu}[2]{{#1}:\ldots:{#2}}

\nc{\maps}[3]{{#1}\!:\!{#2}\rightarrow{#3}}
\nc{\map}[2]{{#1}\rightarrow {#2}} \nc{\res}[2]{{#1} |_{#2}}
\nc{\imbed}{\hookrightarrow}
\title{Kleinian groups of small Hausdorff dimension are classical Schottky groups. I}
\author{ Yong Hou }
\date{}
\date{\small{dedicated to Ying Zhou}}
\begin{document}
\maketitle
\begin{abstract}
It has been conjectured that the Hausdorff dimensions of nonclassical Schottky groups are strictly bounded from below.
In this first part of our works on this conjecture, we prove that there exists a universal positive number $\lambda>0$,  
such that any $2$-generated non-elementary Kleinian groups with limit set of Hausdorff dimension $<\lambda$
are classical Schottky groups. 
\end{abstract}
\setcounter{tocdepth}{1}
\tableofcontents
\section{Introduction and Main Theorem}
Let $\mathbb{H}^3$ be the hyperbolic $3$-space. 
A subgroup $\G$ of $\Psl(2,\mathbb{C})=\iso(\mathbb{H}^3)$ is called a Kleinian group
if it is discrete. Let $x\in\mathbb{H}^3$. The orbit of $x$ under action of $\G$ is
denoted by $\G x$. The limit set $\Lambda_\G$ of $\G$ is defined as
$\Lambda_\G=\overline{\G x}\cap\partial\mathbb{H}^3$. By definition, $\Lambda_G$
is the smallest closed $\G$-invariant subset of $\partial\mathbb{H}^3$. The group
$\G$ is said to be \emph{elementary} if $\Lambda_G$ contains at most two points, otherwise $\G$
is said to be \emph{non-elementary}. Note that elementary Kleinian groups are completely classified.
Henceforth when we say "Kleinian group $\G$" we will assume that $\G$ is non-elementary. The group $\G$
is called second kind if $\Lambda_\G\not=\partial\mathbb{H}^3$, otherwise it
is said to be first kind. The set
$\Omega_\G=\partial\mathbb{H}^3-\Lambda_\G$ is the region of discontinuity, and
$\G$ acts properly discontinuously on $\Omega_\G$. \par
Let $\{\Delta_1,\Delta'_1,...,\Delta_k,\Delta'_k\}$ be a collection of disjoint
closed Jordan curves in the Riemann sphere $\overline{\mathbb{C}}$ and let $D_i,D'_i$ be the topological disks bounded
by $\Delta_i$ and $\Delta'_i$ respectively. Suppose we have 
$\{\g_i\}^k_1\subset\Psl(2,\mathbb{C})$
such that $\g_i(\Delta_i)=\Delta'_i$ and $\g_i(D_i)\cap D'_i=\emptyset.$ Then the group
$\G$ generated by $\{\g_1,...,\g_k\}$ is a free Kleinian
group of rank $k,$ and $\G$ is called marked Schottky group with marking $\{\g_1,...,\g_k\}$. 
A finitely generated Kleinian group $\G$ is called Schottky group if it is a marked Schottky group for
some markings.
If there exists a generating 
set $\{\g_i,...,\g_k\}$ such that all $\Delta_i,\Delta'_i$ can be taken as circles then it is 
called a marked classical Schottky group with classical markings $\{\g_i,...,\g_k\}$, and 
$\{\g_1,...,\g_k\}$ are called classical generators.
A Schottky group $\G$ is called classical Schottky group if there exists a classical markings for $\G$.\par
For a Schottky group
$\G$, the manifold $\mathbb{H}^3/\G$ is homeomorphic to the interior of a handlebody 
of genus $k$. We denote by $\mathfrak{J}_k$ the set of all rank $k$ Schottky groups,
and let $\mathfrak{J}_{k,o}$ be the set of all rank $k$ \emph{classical} Schottky groups. 
One simple way to topologize $\mf{J}_k$ into a topological space is to identify it with the space of moduli for 
a Riemann surface of genus $k$.\par
It is known that not all Schottky groups are classical Schottky groups, in-fact
the space of classical Schottky groups is not even dense in the space of Schottky
groups \cite{Marden}, also see \cite{Doyle}.\par In \cite{Doyle}, Peter Doyle proved that there exists a universal
upper bound on the Hausdorff dimension of the limit sets of finitely generated classical Schottky groups. 
It was originally Phillips and Sarnak in \cite{phillips}
who proved that there exists a universal upper bound on the Hausdorff dimension of the limit sets of classical Schottky groups of dimension greater than $3$.\par
Let $\mathfrak{D}_\G$ denote the Hausdorff dimension of $\Lambda_\G,$  
the main result is the following.
\pagebreak
\begin{thm}\label{main} 
There exists a universal $\lambda>0$, 
such that any $2$-generated non-elementary Kleinian group $\G$ with $\mathfrak{D}_\G<\lambda$ is a classical Schottky group.
\end{thm}
Note that our result can be viewed as the converse of the result by Doyle \cite{Doyle} and
Phillips and Sarnak \cite{phillips}. The proofs of their theorems relies on a crucial fact that 
$\lambda_0(\mathbb{H}^{n+1}/\G)=\mathfrak{D}_\G(n-\mathfrak{D}_\G)$, for $\mathfrak{D}_\G\ge 1,$ where $\lambda_0(\mathbb{H}^{n+1}/\G)$ is the bottom spectrum of the Laplacian on the hyperbolic manifold $\mathbb{H}^{n+1}/\G.$
But this identity obviously
is useless in our situation.\par
We prove Theorem \ref{main} by using a result of \cite{Hou}, and selections of generators.
The proof is divided into three main steps. \par
To lead up to the proof, we first do some preliminary estimates on the locations of the fixed points of a 
given set of generators of a Schottky group. These estimates give us a sufficient control on how the fixed points
of a set of generators change in terms of the Hausdorff dimension of the limit set of the group. The main ingredient 
of the proofs of these estimates relies on the result of \cite{Hou}, Theorem $1.1$ rewritten in the trace form. 
\par
Next we obtain a set of sufficient conditions for any given sequence of Schottky groups to contain a subsequence
of classical Schottky groups in the unit ball in hyperbolic space. These conditions are stated in the upper-half space 
hyperbolic model. The idea is that, if the radius of isometric circles of a sequence
of generators, decreases sufficiently faster than the reduction of the gaps between any of the fixed points of the sequence
of generators, then this sequence of generators will eventually become classical generators.
We do this first by transforming the generators with the condition that 
the generator with the shortest translation length
is transformed into vertical position passing through the origin with fixed points at north and south poles. And 
then these generators are projected into upper-half space. \par

In the first step of the proof, we consider Schottky subspaces of the Schottky space that consist of Schottky groups
$\G$ for which there exists a
generating set $S_\G$ for $\G$ with the set of fixed points of $S_\G$
on the boundary sphere on the unit $3$-ball hyperbolic space that are mutually bounded away by a 
positive constant. \par
In this step we prove that Theorem \ref{main} holds for these Schottky subspaces. 
This is proved via contradiction.  
Suppose $\G_n$ is a sequence of nonclassical Schottky groups 
in the subspace with Hausdorff dimension of $\Lambda_{\G_n}$ decreasing to $0.$ \par 
The idea is that we first transform these generators of the 
generating sets $S_{\G_n}$  into the standard form with the generator of $S_{\G_n}$ of shortest translation length 
put in vertical position. If no generator of $S_{\G_n}$ is of bounded translation length when $\mathfrak{D}_{\G_n}\to 0$ then
it's easy to see it will lead to a contradiction. On the other hand,
it's a simple corollary of \cite{Hou} that there can exists at most one generator of bounded translation length per $S_{\G_n}$
when $\mathfrak{D}_{\G_n}\to 0.$ If such a generator does exists, then we first make a careful change of generators which will be constructed 
based on estimates that it's fixed-point set will be ``minimally excluded" from isometric spheres of 
the generator of bounded translation length. Working with these generating sets we will show that with appropriate
additional transformations and changes the generators with translation length that is not bounded will always grow sufficiently fast 
to lead to disjoint isometric spheres and is also disjoint from the rest of the isometric spheres of the other generators.\par 
One of the crucial tools in estimating this growth is the strong form of the inequality of \cite{Hou}. As we will show that 
when one of the generators is of bounded translation length then it will force 
the needed growth for the rest of the generators. And with appropriate choice of generators this will lead to a classical
Schottky group.\par
The second step, we prove that on the nonclassical Schottky space there exist universal lower bounds for
$\mathfrak{D}_\G+Z_\G,$ see Section $2$ for notation definitions. Essentially this means that we cannot have simultaneous arbitrary small Hausdorff dimension and
minimal gaps of its fixed-points set on the space of nonclassical Schottky groups. This is proved using generator selection
and the results of step one.
\par
In the last step, we prove that when $\mathfrak{D}_{\G}$ is taken sufficiently small then $\G$ can be taken as Schottky group.
This is based on basic topological arguments and some well known results on Kleinian groups. Although we proved 
(easy standard topological argument) that all finitely generated Kleinian groups with Limit set of sufficiently small
Hausdorff dimension is a Schottky group, but as Dick Canary pointed out that it's simpler in the $2$-generated case based on 
the result of Peter Shalen which assures
that any $2$-generated Kleinian group is either free or cofinite volume.
\par
The paper is organized as follows:
In Section $2$, we define and lists global notations which will be used throughout the paper.
In Section $3,$ a strong form of Theorem 1.1 given in \cite{Hou} will be stated for two generators, which we shall
use for selecting generators, see Corollary \ref{cor-2}. In Section $4,$ for a given sequence of Schottky groups $\G_n$ with bounds on
$Z_{\G_n}$ (see section $2$),
we prove inequalities that will enable us to control fixed points
of a given sequence of generators of $\G_n$ in relation to the Hausdorff dimensions $\mathfrak{D}_{\G_n}$ of $\G_n.$
These will be used in the selection processes. In Section $5$, sufficient conditions for a given pair of generators
to be classical generators is proved which we will use in our generator selection process. In Section $6$, we will
use tools developed in previous sections to form a generator selection process and prove that Schottky groups with small $\mathfrak{D}_{\G_n}$ and bounds
on $Z_{\G_n}$ are classical Schottky groups. In Section $7$, we will prove the theorem that will remove the bound condition on
$Z_{\G_n}.$  Section $8,$ completes the proofs of our main theorem by reducing finitely generated Kleinian groups with small Hausdorff dimensions to
Schottky groups via a standard known topological argument.\par
\centerline{ACKNOWLEDGEMENTS}\par
I wish to express my deepest appreciation to the referee for 
spent enormous amount of time reading and correcting the paper which has made the paper much better.
I would like to express my sincere appreciation and gratitude to Benson Farb for providing me guidance and support in the
write-up of the paper.
I would like to express my gratitude to Jim W. Anderson, Peter Shalen, Marc Culler, Dick Canary for their interests on the paper. 
I am indebted to Ying Zhou for all the support and unwavering help that has been provided to me during this work.\par

\section{Notations}
Let $\G\subset\Psl(2,\mathbb{C})$ denote a Schottky group generated by $<\alpha,\beta>$ with $\alpha$ having fixed points $0,\infty.$ 
Asssume that $\gamma\in\G$ is a loxodromic element having fixed points
$\not=\infty.$ Write $\gamma$ in matrix form 
$\gamma=\bigl(\begin{smallmatrix} a & b\\c & d\end{smallmatrix}\bigr),$ with $\det(\gamma)=1.$
We will set the following notations and definitions throughout the rest of the paper.

\begin{nota}\text{}
\begin{itemize}
\item Denote the critical exponent of $\G$ by $D_\G,$ and the Hausdorff dimension of $\Lambda_\G$ by $\mathfrak{D}_\G.$ 
\item $\mathfrak{R}_\gamma$ the radius of isometric circles of $\gamma$.
\item $\eta_\gamma=-\frac{d}{c}$, $\zeta_\gamma=\frac{a}{c}.$
\item We will define two different ways to denotes the two fixed points of $\gamma:$ 
$\left\{z_{\gamma,l},z_{\gamma,u}\right\}$ the two fixed
points of $\gamma$ in $\mathbb{C}$ with $|z_{\gamma,l}|\le|z_{\gamma,u}|$, and $\left\{z_{\gamma,-},z_{\gamma,+}\right\}$ 
the two fixed points with $z_{\gamma,\pm}$ given by quadratic formula with subscripts $\pm$
corresponding to $\pm\sqrt{\tr^2(\gamma)-4}.$ Note that we always take the principle branch for the square roots of complex numbers.
\item $\mathcal{L}_\gamma$ is the axis of $\gamma.$
\item $T_\gamma$ is the translation length of $\gamma.$
\item Set $Z_{\beta}:=\min\{|z_{\beta,-}-z_{\beta,+}|,|\frac{1}{z_{\beta,-}}-\frac{1}{z_{\beta,+}}|\}.$ 
\item
$Z_{<\alpha,\beta>}:=\min\{Z_{\beta},|z_{\beta,+}|,|z_{\beta,-}|,|z_{\beta,+}|^{-1},|z_{\beta,-}|^{-1}\}.$
\item For $\epsilon>0,$ we say that $Z_\G>\epsilon,$ if there exists a generating set $<\alpha,\beta>$ of $\G$ such that
$Z_{<\alpha,\beta>}>\epsilon.$
\item
Given any two sequences of real numbers $\{p_n,q_n\}$, we use the notation $p_n\asymp q_n$ iff, 
there exists $\sigma>0$ such that, $\sigma^{-1}<\liminf\frac{p_n}{q_n}\le\limsup\frac{p_n}{q_n}<\sigma.$\par

\end{itemize}
\end{nota}
\begin{nota}
Let $\{\gamma_n\}\subset\Psl(2,\mathbb{C})$ be a sequence of loxodromic transformations.
Let $\{p_n\}$ be a sequence of complex numbers, and $\{q_n\}$ a sequence of positive real numbers. We write:
\[\left|z_{\gamma_n,\pm}-p_n\right|<q_n\]
if there exists $N$ such that for every $n>N$ we have at least one of the following holds,
\begin{itemize}
\item[(i)]
\[\left|z_{\gamma_n,+}-p_n\right|<q_n,\]
\item[(ii)]
\[\left|z_{\gamma_n,-}-p_n\right|<q_n.\]
\end{itemize}
\end{nota}

\section{Free Group Actions}
Given a finitely generated non-elementary Kleinian group $\G$, the critical exponent of $\G$ is the unique positive 
number $D_\G$ such that the Poincar\'{e} series of $\G$ given by 
$\sum_{\g\in\G}e^{-s\dis(x,\g x)}$ is divergent if $s<D_\G$ and
convergent if $s>D_\G$. If the Poincar\'{e} series diverges at
$s=D_\G$ then $\G$ is said to be divergent. Bishop-Jones showed that
$D_\G\le\mathfrak{D}_\G$ for all analytically finite non-elementary Kleinian groups $\G$.
In-fact, if $\G$ is topologically tame ($\mathbb{H}^3/\G$ homeomorphic to the
interior of a compact manifold-with-boundary) then  $D_\G=\mathfrak{D}_\G$. Hence it follows
Agol's \cite{agol} proof of tameness conjecture that $D_\G=\mathfrak{D}_\G$ for all finitely generated non-elementary Kleinian groups.
The critical exponent is a geometrically rigid object in the sense that
a decrease in $D_\G$ corresponds to a decrease in geometric complexity.\par
Next we state the following theorem from \cite{Hou}, which provides the
relation between the group action and the critical exponent.
\begin{thm}[Hou\cite{Hou}]\label{critical}
Let $\G$ be a free non-elementary Kleinian group of rank $k$ with free generating
set $\mathcal{S}$, and $x\in\mathbb{H}^3$ then 
\[ \sum_{\g\in\mathcal{S}}\frac{1}{1+\exp(D_\G\dis(x,\g x))}\le\frac{1}{2}.\]
In particular we have at least $k-1$ distinct elements $\{\g_{i_j}\}_{1\le j\le k-1}$ 
of $\mathcal{S}$ that satisfies $\dis(x,\g_{i_j} x)\ge\log(3)/D_\G,$ and at least one element
$\g_{i_j}$ with $\dis(x,\g_{i_j} x)\ge\log(2k-1)/D_\G.$
\end{thm}
The following is a useful corollary of Theorem \ref{critical}, stated
here for the case of $\G$ is a free group of rank $2$.
\begin{cor}\label{cor-2}
Let $\mathcal{S}=\{\g_1,\g_2\}$ be a generating set for a free non-elementary Kleinian group $\G$. Let
$x\in\mathbb{H}^3$. Then
\[ \dis(x,\g_2 x)\ge\frac{1}{D_{\G}}\log\left(\frac{e^{D_\G\dis(x,\g_1 x)}+3}{e^{D_\G\dis(x,\g_1 x)}-1}\right)\]
\end{cor}
\begin{cor}\label{cor-1}
Let $\mathcal{S}=\{\g_1,\g_2\}$ be a generating set for a free non-elementary Kleinian group $\G$. Let
$x\in\mathbb{H}^3$.
Let $m$ be any integer. Then at least one of the elements $\g'$ of 
$\mathcal{S}'=\{\g^m_1\g_2,\g^{m+1}_1\g_2\}$ 
satisfies $\dis(x,\g' x)\ge\log 3/D_\G$.
\end{cor}
\section{Trace, Fixed Points and Hausdorff Dimension}
In this section we study the relationships of fixed points of generating sets of a given sequence of Schottky groups $\G_n$
and the Hausdorff dimensions of $\Lambda_{\G_n}.$\par
What we like to do is
to find a relationship between the distribution of the fixed points of one of the generators in terms
of the translation length growth of the other generator and the Hausdorff dimension of its limit set. By having
this type of relationship we will be able to construct a new set of generators from the
given generating set with prescribed distribution of its fixed points. The new set of
generators will be a crucial ingredient in the proof of our theorem.\par
Let $\{\G_n\}$ be a sequence of rank $2$ Schottky groups with $D_n\to 0$ generated by $\alpha_n,\beta_n\in\Psl(2,\mathbb{C})$
in the upper space model $\mathbb{H}^3$ with 
$\alpha_n=\bigl(\begin{smallmatrix} \lambda_n & 0\\0 & \lambda^{-1}_n\end{smallmatrix}\bigr)$, $|\lambda_n|>1$.
Set $\beta_n=\bigl(\begin{smallmatrix} a_n & b_n\\c_n & d_n\end{smallmatrix}\bigr).$ \par
We assume throughout this section that there exists $M>0$ such that $T_{\alpha_n}<M$ for all $n$. Set $D_n=D_{\G_n}.$ Let $\tr=\text{trace}.$\par
First we will state Corollary \ref{cor-2} in the trace form. 
\begin{prop}\label{trace}
Suppose there exists $\Delta>0$ such that $Z_{\beta_n}>\Delta.$
There exists $\rho>0$ depending on $\Delta$, such that
\[|\tr(\beta_n)|>\rho\left(\frac{|\lambda_n|^{2D_n}+3}{|\lambda_n|^{2D_n}-1}\right)^{\frac{1}{2D_n}}\]
for large $n$.
\end{prop}
\begin{proof}
Let $T_n$ be the translation length of $\beta_n$ and $R_n=\dis(\mathcal{L}_{\alpha_n},\mathcal{L}_{\beta_n}).$ 
Let $x_n$ be a point on axis of $\alpha_n$ which is the nearest point of $\mathcal{L}_{\alpha_n}$ to $\mathcal{L}_{\beta_n}.$ 
By triangle inequality, 
$T_n\ge \dis(x_n,\beta_nx_n)-2R_n$. And for sufficiently large $T_n$, we have for some positive constant
$c>0,$ $|\tr^2(\beta_n)|>ce^{T_n}.$ Now for large $n,$ from Corollary \ref{cor-2},
\[|\tr^2(\beta_n)|>c\left(\frac{|\lambda_n|^{2D_n}+3}{|\lambda_n|^{2D_n}-1}\right)^{\frac{1}{D_n}}
\left(e^{-2\dis(\mathcal{L}_{\alpha_n},\mathcal{L}_{\beta_n})}\right). \] 
Now $Z_{\beta_n}>\Delta$ implies that $\dis(\mathcal{L}_{\alpha_n},\mathcal{L}_{\beta_n})<M$ for some $M>0.$ Hence the result follows.
\end{proof}
\begin{remb}\label{4A}
Note that without assuming bounds on $Z_{\beta_n}$ we can state above Proposition \ref{trace} as,
\[|\tr^2(\beta_n)|>c\left(\frac{|\lambda_n|^{2D_n}+3}{|\lambda_n|^{2D_n}-1}\right)^{\frac{1}{D_n}}
\left(e^{-2\dis(\mathcal{L}_{\alpha_n},\mathcal{L}_{\beta_n})}\right). \] 
If $\dis(\mathcal{L}_{\alpha_n},\mathcal{L}_{\alpha_n\beta_n})<\epsilon$ then there exists $\delta>0$ such that,
\[|\tr(\beta_n)|>\rho\left(\frac{|\lambda_n|^{2D_n}+3}{|\lambda_n|^{2D_n}-1}\right)^{\frac{1}{2D_n}}\]
for large $n$.
\end{remb}
The next lemma and its corollaries are estimates of convergence rates of fixed points of the generators of $\G_n$  in terms of the 
Hausdroff dimension of $\Lambda_n.$
\begin{lem}\label{fix-trace}
Suppose there exists $\Delta>0,M>0$ such that $Z_{<\alpha_n,\beta_n>}>\Delta$ and $T_{\alpha_n}<M$ for all $n.$ 
Let $k_n,l_n$ be any integers such that $T_{{\alpha_n}^{k_n}},T_{{\alpha_n}^{l_n}}<M$. Then
there exists a constant $\rho>0$ such that, 
\[\left|z_{\alpha^{k_n}_n\beta_n\alpha^{l_n}_n,\pm}-\zeta_{\alpha^{k_n}_n\beta_n\alpha^{l_n}_n}\right|+
\left|z_{\alpha^{k_n}_n\beta_n\alpha^{l_n}_n,\mp}-\eta_{\alpha^{k_n}_n\beta_n\alpha^{l_n}_n}\right|
<\frac{\rho}{|\tr(\beta_n)||\lambda^{l_n-k_n}_n|}\]
for large $n.$
\end{lem}
\begin{rem}\label{rem-1}
Lemma \ref{fix-trace} is an estimate of how fast the fixed points converges, it's not important to our applications in this paper
which fixed point converges to $\eta_\g$ and which converges to $\zeta_\g$ for a given $\g\in\Psl(2,\mathbb{C}).$ \par 
Given a complex number $z=re^{i\theta}$ we write $\sqrt{z^2}=z$ if $-\pi<2\theta\le\pi,$ and $\sqrt{z^2}=-z$ if $2\theta>\pi$ or $2\theta\le-\pi.$
Then a more precise statement of Lemma \ref{fix-trace}
which dichotomizes the above inequality for large $n$ would be:
\begin{itemize}
\item[(i)] $\sqrt{\tr^2(\alpha^{k_n}_n\beta_n\alpha^{l_n}_n)}=\tr(\alpha^{k_n}_n\beta_n\alpha^{l_n}_n)$ 
\[\left|z_{\alpha^{k_n}_n\beta_n\alpha^{l_n}_n,+}-\zeta_{\alpha^{k_n}_n\beta_n\alpha^{l_n}_n}\right|+
\left|z_{\alpha^{k_n}_n\beta_n\alpha^{l_n}_n,-}-\eta_{\alpha^{k_n}_n\beta_n\alpha^{l_n}_n}\right|
<\frac{\rho}{|\tr(\beta_n)||\lambda^{l_n-k_n}_n|}\]

\item[(ii)] $-\sqrt{\tr^2(\alpha^{k_n}_n\beta_n\alpha^{l_n}_n)}=\tr(\alpha^{k_n}_n\beta_n\alpha^{l_n}_n)$
\[\left|z_{\alpha^{k_n}_n\beta_n\alpha^{l_n}_n,-}-\zeta_{\alpha^{k_n}_n\beta_n\alpha^{l_n}_n}\right|+
\left|z_{\alpha^{k_n}_n\beta_n\alpha^{l_n}_n,+}-\eta_{\alpha^{k_n}_n\beta_n\alpha^{l_n}_n}\right|
<\frac{\rho}{|\tr(\beta_n)||\lambda^{l_n-k_n}_n|}\]

\end{itemize}
\end{rem}
\begin{cor}\label{fix-bound-1}
Suppose there exists $\Delta>0,M>0$ such that $Z_{<\alpha_n,\beta_n>}>\Delta$ and $T_{\alpha_n}<M$ for all $n.$ 
Let $k_n,l_n$ be any integers such that $T_{{\alpha_n}^{k_n}},T_{{\alpha_n}^{l_n}}<M$. Then
for any $\delta>0$ there exists $\epsilon>0$ such that if $D_n<\epsilon$ then,
\[\left|z_{\alpha^{k_n}_n\beta_n\alpha^{l_n}_n,\pm}-\zeta_{\alpha^{k_n}_n\beta_n\alpha^{l_n}_n}\right|+
\left|z_{\alpha^{k_n}_n\beta_n\alpha^{l_n}_n,\mp}-\eta_{\alpha^{k_n}_n\beta_n\alpha^{l_n}_n}\right|<\delta(|\lambda_n|^2-1).\]
The same dichotomy decomposition of the inequality holds as given in Remark \ref{rem-1}.
\end{cor}
\begin{proof}
Let us assume that $\sqrt{\tr^2(\alpha^{k_n}_n\beta_n\alpha^{l_n}_n)}=\tr(\alpha^{k_n}_n\beta_n\alpha^{l_n}_n).$
Also $Z_{<\alpha_n,\beta_n>}\ge\Delta$ by Proposition \ref{trace} as $D_n\to 0$ we have $|\tr(\beta_n)|\to\infty.$ In addition, by 
$Z_{<\alpha_n,\beta_n>}\ge\Delta,$ we have $|z_{+,\beta_n}-z_{-,\beta_n}|\ge\Delta,$ and also
$|z_{+,\beta_n}|^{-1}, |z_{-,\beta_n}|^{-1}\ge\Delta$ gives $|z_{+,\beta_n}-z_{-,\beta_n}|\le|z_{+,\beta_n}|+|z_{-,\beta_n}|<2/\Delta.$
There exists $c_1,c_2>0$
with $c_1<|z_{+,\beta_n}-z_{-,\beta_n}|<c_2.$ By $|z_{+,\beta_n}-z_{-,\beta_n}|=|\frac{\sqrt{\tr^2(\beta_n)-4}}{2c_n}|,$ 
such we have $c_1<|\frac{\sqrt{\tr^2(\beta_n)-4}}{2c_n}|<c_2.$ This implies $|\tr(\beta_n)|\asymp |c_n|,$
and since $\tr(\beta_n)=a_n+d_n$, there exists $c_3,c_4>0$ such that $c_3<|(a_n+d_n)/c_n|<c_4.$ Since $|\tr(\beta_n)|\to\infty,$ 
implies $\mathfrak{R}_{\beta_n}\to 0,$ and $Z_{<\alpha_n,\beta_n>}\ge\Delta$  gives us that
there exists $c_5,c_6>0$ such that $c_5<|\frac{a_n}{c_n}|,|\frac{d_n}{c_n}|<c_6,$ for large $n.$ \par
Therefore the fixed points $z_{\pm,\beta_n}$
of $\beta_n$ must $\to\{\frac{a_n}{c_n},-\frac{d_n}{c_n}\}$. There exists $c,c',N>0$ such
that:
\[\left|z_{\alpha^{k_n}_n\beta_n\alpha^{l_n}_n,+}-\frac{a_n}{c_n}\lambda^{2k_n}_n\right|=
\left|\frac{(a_n\lambda^{k_n+l_n}_n-d_n\lambda^{-l_n-k_n}_n)+\sqrt{\tr^2(\alpha^{k_n}_n\beta_n\alpha^{l_n}_n)-4}}
{2c_n\lambda^{l_n-k_n}_n}-\frac{a_n\lambda^{k_n+l_n}_n}{c_n\lambda^{l_n-k_n}_n}\right|\]
\[=\left|\frac{\sqrt{\tr^2(\alpha^{k_n}_n\beta_n\alpha^{l_n}_n)-4}-\sqrt{\tr^2(\alpha^{k_n}_n\beta_n\alpha^{l_n}_n)}}
{2c_n\lambda^{l_n-k_n}_n}\right|\]
\[\le \frac{c}{\left|\lambda^{l_n-k_n}_n\right||\tr(\beta_n)|}\left|\sqrt{\tr^2(\alpha^{k_n}_n\beta_n\alpha^{l_n}_n)-4}-
\sqrt{\tr^2(\alpha^{k_n}_n\beta_n\alpha^{l_n}_n)}\right|\]
\[\text{if}\quad |\tr(\alpha^{k_n}_n\beta_n\alpha^{l_n}_n)|<\kappa,\quad \text{for some, }\kappa>0, \text{ and all}\quad  n\quad \text{then,}\]
\[\le\frac{c\kappa'}{\left|\lambda^{l_n-k_n}_n\right||\tr(\beta_n)|} \quad\text{for some $\rho'>0$.}\]
Otherwise we write,\hspace{8cm}
\[\le c\frac{|\sqrt{\tr^2(\alpha^{k_n}_n\beta_n\alpha^{l_n}_n)-4}|}{\left|\lambda^{l_n-k_n}_n\right||\tr(\beta_n)|}\left|\frac{1}{\sqrt{1-4/\tr^2(\alpha^{k_n}_n\beta_n\alpha^{l_n}_n)}}-1\right|\]
by using Binomial series,\hspace{8cm}
\[(1-4/\tr^2(\alpha^{k_n}_n\beta_n\alpha^{l_n}_n))^{-1/2}=1+2\tr^{-2}(\alpha^{k_n}_n\beta_n\alpha^{l_n}_n)+\epsilon_n,\]
with $\epsilon_n\to 0$ at order $|\tr^{-4}(\alpha^{k_n}_n\beta_n\alpha^{l_n}_n)|$. Hence there $\sigma>0$ such that,
\[\le\frac{c\sigma}{\left|\lambda^{l_n-k_n}_n\right|\left|\tr(\beta_n)\right|\left|\tr(\alpha^{k_n}_n\beta_n\alpha^{l_n}_n)\right|}\quad\quad\text{for large $n$.}\]
Hence in either case we have for some $c'>0$ such that,\hspace{3cm}
\[\le\frac{c'}{\left|\lambda^{l_n-k_n}_n\right|\left|\tr(\beta_n)\right|}\quad\quad\text{for}\quad n>N.\]
This gives (i) of the Lemma. \par
Now $T_{{\alpha_n}^{k_n}},T_{{\alpha_n}^{l_n}}<M,$ implies $\left|\lambda^{l_n-k_n}_n\right|<M'$ for $M'>0.$ Hence,
\[\left|z_{\alpha^{k_n}_n\beta_n\alpha^{l_n}_n,+}-\frac{a_n}{c_n}\lambda^{2k_n}_n\right|\le\frac{c''}{|\tr(\beta_n)|},\quad n>N, c''>0.\]
By $4\cosh(T_{\beta_n})=|\tr^2(\beta_n)|+|\tr^2(\beta_n)-4|,$ we have
\[e^{T_{\beta_n}}<|\tr^2(\beta_n)|(1+2/|\tr^2(\beta_n)|).\]
Since $|\tr(\beta_n)|\to\infty,$ there exists $a>0$ such that $e^{T_{\beta_n}}\le a|\tr(\beta_n)|^2$ for sufficient large $n.$ \par
Note that if $|\tr(\alpha^{k_n}_n\beta_n\alpha^{l_n}_n)|\to\infty$ then we also have the stronger inequality,
\[\left|z_{\alpha^{k_n}_n\beta_n\alpha^{l_n}_n,+}-\frac{a_n}{c_n}\lambda^{2k_n}_n\right|\le\frac{c'''}{|\tr(\beta_n)|
|\tr(\alpha^{k_n}_n\beta_n\alpha^{l_n}_n)|},\quad n>N, c'''>0.\]
Here we can take $c'''>\frac{c\sigma}{|\lambda_n^{l_n-k_n}|}$ as given in the above binomial inequality.\\ 
And $b>0$ such that $e^{T_{\alpha^{k_n}_n\beta_n\alpha^{l_n}_n}}\le b|\tr(\alpha^{k_n}_n\beta_n\alpha^{l_n}_n)|^2.$ We 
therefore need 
\[\min\{\rho e^{-\frac{T_{\beta_n}}{2}},\rho e^{-\frac{T_{\alpha^{k_n}_n\beta_n\alpha^{l_n}_n}}{2}}\}\\<\delta(|\lambda_n|^2-1),\quad\text{for large } n,\text{ and some }\rho>0.\] 
By $|\lambda_n|^2=e^{T_{\alpha_n}}$, we require 
\begin{equation*}
\max\{T_{\beta_n},T_{\alpha^{k_n}_n\beta_n\alpha^{l_n}_n}\}\ge 2\log\left(\frac{\frac{\rho}{\delta}}{e^{T_{\alpha_n}}-1}\right).
\end{equation*}
Hence it follows 
there exists $\rho'>0$ such that if $D_n<\rho'$ then at least one of $T_{\beta_n}, T_{\alpha^{k_n}_n\beta_n\alpha^{l_n}_n}$
satisfies the above inequality. Therefore we have $|z_{\beta_n,+}-\frac{a_n}{c_n}|<\delta(|\lambda_n|^2-1)$
for $D_n<\rho$. The proof for the other part is same.\par 
The case $-\sqrt{\tr^2(\alpha^{k_n}_n\beta_n\alpha^{l_n}_n)}=\tr(\alpha^{k_n}_n\beta_n\alpha^{l_n}_n)$ is similar with
$z_{\alpha^{k_n}_n\beta_n\alpha^{l_n}_n,+}$ replace with $z_{\alpha^{k_n}_n\beta_n\alpha^{l_n}_n,-}$ and vice versa.
\end{proof}
\begin{rem}\label{rem-2}
Based on the proof above, we note that the condition $Z_{<\alpha_n,\beta_n>}>\Delta$ in Lemma \ref{fix-trace} can be replaced
with conditions $|\tr(\beta_n)|\to\infty$ and $|\tr(\beta_n)|\le|c_n|.$
\end{rem}
\begin{cor}\label{2-tr-fixed}
Suppose there exists $\Delta>0,M>0$ such that $Z_{<\alpha_n,\beta_n>}>\Delta$ and $T_{\alpha_n}<M$ for all $n.$ 
Let $k_n,l_n$ to be any integers such that $T_{{\alpha_n}^{k_n}},T_{{\alpha_n}^{l_n}}<M$.
Then there exists constants $\delta, \rho>0$ such that, if $|\tr(\alpha^{k_n}_n\beta_n\alpha^{l_n}_n)|\to\infty$
then,
\[\left|z_{\alpha^{k_n}_n\beta_n\alpha^{l_n}_n,\pm}-\zeta_{\alpha^{k_n}_n\beta_n\alpha^{l_n}_n}\right|+
\left|z_{\alpha^{k_n}_n\beta_n\alpha^{l_n}_n,\mp}-\eta_{\alpha^{k_n}_n\beta_n\alpha^{l_n}_n}\right|
<\frac{\rho}{|\tr(\alpha^{k_n}_n\beta_n\alpha^{l_n}_n)||\tr(\beta_n)|}\]
for all $D_n<\delta$.
The same dichotomy decomposition of the inequality holds as given in Remark \ref{rem-1}.
\end{cor}
\begin{proof}
Let us assume that $\sqrt{\tr^2(\alpha^{k_n}_n\beta_n\alpha^{l_n}_n)}=\tr(\alpha^{k_n}_n\beta_n\alpha^{l_n}_n).$ 
Using the inequality in the proof of Corollary \ref{fix-bound-1} and we take 
$\sigma'>\frac{c}{|\lambda_n^{l_n-k_n}|},$ with $c$ given in the proof of Corollary \ref{fix-bound-1},
\[
\left|z_{\alpha^{k_n}_n\beta_n\alpha^{l_n}_n,+}-\frac{a_n}{c_n}\lambda^{2k_n}_n\right|\le
\sigma'\frac{|\sqrt{\tr^2(\alpha^{k_n}_n\beta_n\alpha^{l_n}_n)-4}|}{|\tr(\beta_n)|}\left|\frac{1}
{\sqrt{1-4/\tr^2(\alpha^{k_n}_n\beta_n\alpha^{l_n}_n)}}-1\right|\]
by using the binomial series as in the proof of Corollary \ref{fix-bound-1} and take $\sigma''>0$ as 
$\sigma''>\frac{\sigma\sigma'}{|\lambda_n^{l_n-k_n}|}$ here $\sigma$ is the constant given in the proof of Corollary \ref{fix-bound-1} earlier we have,
\[\left|z_{\alpha^{k_n}_n\beta_n\alpha^{l_n}_n,+}-\frac{a_n}{c_n}\lambda^{2k_n}_n\right|\le\frac{\sigma''}{|\tr(\alpha^{k_n}_n\beta_n\alpha^{l_n}_n)||\tr(\beta_n)|}.\]
The proof for the other part is same.\par
The case $-\sqrt{\tr^2(\alpha^{k_n}_n\beta_n\alpha^{l_n}_n)}=\tr(\alpha^{k_n}_n\beta_n\alpha^{l_n}_n)$ is similar with
$z_{\alpha^{k_n}_n\beta_n\alpha^{l_n}_n,+}$ replaced with $z_{\alpha^{k_n}_n\beta_n\alpha^{l_n}_n,-}$ and vice versa.
\end{proof}
\begin{lem}\label{fix-bound-2}
Suppose there exists $\Delta>0,M>0$ such that $Z_{<\alpha_n,\beta_n>}>\Delta$ and $T_{\alpha_n}<M$ for all $n.$ 
Let $k_n,l_n$ to be any integers such that $T_{{\alpha_n}^{k_n}},T_{{\alpha_n}^{l_n}}<M$.
Then there
exists a constant $\sigma_1,\sigma_2>0$ such that,
\[\frac{\sigma_1|\tr(\alpha^{k_n}_n\beta_n\alpha^{l_n}_n)|}{|\tr(\beta_n)|}\ge
|z_{\alpha^{k_n}_n\beta_n\alpha^{l_n}_n,+}-z_{\alpha^{k_n}_n\beta_n\alpha^{l_n}_n,-}|
\ge\frac{\sigma_2|\tr(\alpha^{k_n}_n\beta_n\alpha^{l_n}_n)|}{|\tr(\beta_n)|},\]
for all $n$ sufficiently large.
\end{lem}
\begin{proof}
Note that we have,
\begin{eqnarray*}
\left|\frac{a_n}{c_n}\lambda^{2k_n}_n-(-\frac{d_n}{c_n}\lambda^{-2l_n}_n)\right|&=&
\left|\frac{a_n\lambda^{k_n+l_n}_n+d_n\lambda^{-k_n-l_n}_n}{c_n\lambda^{l_n-k_n}_n}\right|
=\left|\frac{\tr(\alpha^{k_n}_n\beta_n\alpha^{l_n}_n)}{c_n\lambda^{l_n-k_n}_n}\right|\label{equal}
\end{eqnarray*}
As earlier in the proof of Lemma \ref{fix-trace} we have, $Z_{<\alpha_n,\beta_n>}>\Delta$ implies $c<\frac{\sqrt{\tr^2(\beta_n)-4}}{2c_n}<c'$ 
for some $c,c'>0$ and
$|\tr(\beta_n)|\to\infty$ implies $|\tr(\beta_n)|\asymp c_n$ (i.e. $\kappa c_n<|\tr(\beta_n)|<\kappa' c_n$ for some $\kappa, \kappa'>0$).
Since $|\lambda^{l_n-k_n}|<M'$ we have,
\[\frac{\kappa|\tr(\alpha^{k_n}_n\beta_n\alpha^{l_n}_n)|}{M'|\tr(\beta_n)|}<
\left|\zeta_{\alpha^{k_n}_n\beta_n\alpha^{l_n}_n}-\eta_{\alpha^{k_n}_n\beta_n\alpha^{l_n}_n}\right|<\kappa'\frac{|\tr(\alpha^{k_n}_n\beta_n\alpha^{l_n}_n)|}
{|\tr(\beta_n)|}.\]
By Lemma \ref{fix-trace} the result follows.
\end{proof}
Although next Lemma is not used in the rest of the paper but we include it here to demonstrate relations between fixed points and Hausdorff dimensions.
\begin{lem}\label{fix-bound-3}
Suppose there exists $\Delta>0,M>0$ such that $Z_{<\alpha_n,\beta_n>}>\Delta$ and $M^{-1}<T_{\alpha_n}<M$ for all $n.$ 
Let $k_n,l_n$ to be any integers such that $T_{{\alpha_n}^{k_n}},T_{{\alpha_n}^{l_n}}<M$.
Then for at least one $i\in\{0,1\}$, 
and any integers $k'_n, l'_n $ with $|(k_n-k'_n)|+|(l_n-l'_n)|=i$, we have
$|z_{\alpha^{k'_n}_n\beta_n\alpha^{l'_n}_n,+}-z_{\alpha^{k'_n}_n\beta_n\alpha^{l'_n}_n,-}|
\ge\frac{\kappa}{D_n|\tr(\beta_n)|}$,
for all $n$ sufficiently large. 
\end{lem}
\begin{proof}
Since $T_{\alpha_n}>M^{-1}$ we have $\lambda_n\not\to 1.$ If $|\tr(\alpha^{k_n}_n\beta_n\alpha^{l_n}_n)|<M$ then
$|\zeta_{\alpha^{k_n}_n\beta_n\alpha^{l_n}_n}-\eta_{\alpha^{k_n}_n\beta_n\alpha^{l_n}_n}|\to 0,$ otherwise
we have $Z_{<\alpha_n,\alpha^{k_n}_n\beta_n\alpha^{l_n}_n>}>\Delta$ for some $\Delta>0$ which implies by Proposition \ref{trace},
$|\tr(\alpha^{k_n}_n\beta_n\alpha^{l_n}_n)|\to\infty.$ Hence we have 
$|\lambda_n^2\zeta_{\alpha^{k_n}_n\beta_n\alpha^{l_n}_n}-\eta_{\alpha^{k_n}_n\beta_n\alpha^{l_n}_n}|\not\to 0.$
Since $\eta_{\alpha^{k_n+1}_n\beta_n\alpha^{l_n}_n}=\eta_{\alpha^{k_n}_n\beta_n\alpha^{l_n}_n}$ and 
$\zeta_{\alpha^{k_n+1}_n\beta_n\alpha^{l_n}_n}=\lambda_n^2\zeta_{\alpha^{k_n}_n\beta_n\alpha^{l_n}_n},$
by Lemma \ref{fix-trace} we get $|z_{\alpha^{k_n+1}_n\beta_n\alpha^{l_n}_n,+}-z_{\alpha^{k_n+1}_n\beta_n\alpha^{l_n}_n,-}|>\kappa.$
Then by Lemma \ref{fix-bound-2}, $|\tr(\alpha^{k_n+1}_n\beta_n\alpha^{l_n}_n)|>\kappa'|\tr(\beta_n)|.$ Since $|\tr(\beta_n)|>\log 3/D_n$
we have $|\tr(\alpha^{k_n+1}_n\beta_n\alpha^{l_n}_n)|>\kappa'\log 3/D_n.$ Similarly we also have
$|\tr(\alpha^{k_n}_n\beta_n\alpha^{l_n+1}_n)|>\kappa''\log 3/D_n.$ The result follows from Lemma \ref{fix-bound-2}.\par
We note that the condition $T_{\alpha_n}>M^{-1}$ is of convenience only not necessary, the lemma still holds without this condition.
\end{proof}
Finally in this section we make the following observation based the proof of Lemma \ref{fix-trace}.
Note that if we don't care about the precise upper bound of fixed points of $\alpha^{k_n}_n\beta_n\alpha^{l_n}_n$ to 
$\zeta_{\alpha^{k_n}_n\beta_n\alpha^{l_n}_n},\eta_{\alpha^{k_n}_n\beta_n\alpha^{l_n}_n},$ then we can relax the conditions in Lemma \ref{fix-trace}
and state as follows,
\begin{remc}\label{4B}
Suppose that $|z_{\beta_n,+}-z_{\beta_n,-}|<c$ and $|c_n\lambda_n^{l_n-k_n}|\to\infty.$ Then
$\{z_{\alpha^{k_n}_n\beta_n\alpha^{l_n}_n,+},z_{\alpha_n^{k_n}\beta_n\alpha_n^{l_n},-}\} \to
\{\zeta_{\alpha^{k_n}_n\beta_n\alpha^{l_n}_n},\eta_{\alpha_n^{k_n}\beta_n\alpha_n^{l_n}}\}.$
\end{remc}
\section{Sufficient Conditions}
In this section we will state and prove a set of conditions for a given sequence of Schottky groups with decreasing Hausdorff dimensions
that will be sufficient for the sequence to contain a subsequence of classical Schottky groups.\par
Let $(\mathbb{B}, \dis_B)$ be the unit ball model of hyperbolic $3-$space. 
Let $\pi:\mathbb{B}\longrightarrow \mathbb{H}^3$ be the stereographic hyperbolic isometry. \par
Given a loxodromic element $\alpha,$ of $\Psl(2,\mathbb{C})$ acting on the unit ball $\mathbb{B}$ model of hyperbolic $3$-space,
denote by $S_{\alpha,r}$ and $S_{\alpha^{-1},r}$ the isometric spheres of
Euclidean radius $r$ of $\alpha.$ We set $\lambda_{\pi_*(\alpha)}$ as the multiplier of $\pi_*(\alpha)$ in 
the upper space model $\mathbb{H}^3.$
For $R>0,$ set $C_R$ as the circle in $\mathbb{C}$ about origin of radius $R.$
\begin{prop}\label{proj}
Let $\alpha$ be a loxodromic element of $\Psl(2,\mathbb{C})$ acting on $\mathbb{B},$
with axis passing through the origin and fixed points on north and south poles. Then
$\pi(S_{\alpha,r}\cap\partial\mathbb{B}),\pi(S_{\alpha^{-1},r}\cap\partial\mathbb{B})$ maps to
$C_{\frac{1}{\lambda_{\pi_*(\alpha)}}},C_{\lambda_{\pi_*(\alpha)}}.$
\end{prop}
\begin{proof}
Let $T_\alpha$ be the translation length of $\alpha$. The Euclidean radius $r$ is given by
$r^{-1}=\sinh(\frac{T_\alpha}{2}).$ In terms of $\lambda_{\pi_*(\alpha)},$
\[ r=\frac{2|\lambda_{\pi_*(\alpha)}|}{|\lambda_{\pi_*(\alpha)}|^2-1}\]
Set ${\bf e}=(0,0,1)$ as the north pole of $\partial\mathbb{B}.$
Let $\delta_\alpha,\delta_{\alpha^{-1}}>0$ denote the radius of $\pi(S_{\alpha,r}\cap\partial\mathbb{B})$ and
$\pi(S_{\alpha^{-1},r}\cap\partial\mathbb{B})$ respectively. Then for $x\in S_{\alpha,r}\cap\partial\mathbb{B},$ 
$\delta_\alpha$ is given by,
\begin{eqnarray*}
\delta^2_\alpha=\frac{4r^2}{(1+r^2)|x-{\bf e}|^4}\\
\text{and}\quad\quad |x-{\bf e}|^2=\frac{4}{|\lambda_{\pi_*(\alpha)}|^2+1}\\
\text{this implies}\quad \quad\delta_{\alpha}=|\lambda_{\pi_*(\alpha)}|.
\end{eqnarray*}
Similarly for $x\in S_{\alpha^{-1},r}\cap\partial\mathbb{B}$ we have,
\begin{eqnarray*}
\delta^2_{\alpha^{-1}}=\frac{4r^2}{(1+r^2)|x-{\bf e}|^4}\\
\text{and}\quad\quad |x-{\bf e}|^2=\frac{4|\lambda_{\pi_*(\alpha)}|^2}{|\lambda_{\pi_*(\alpha)}|^2+1}\\
\text{which gives}\quad \quad\delta_{\alpha^{-1}}=\frac{1}{|\lambda_{\pi_*(\alpha)}|}.
\end{eqnarray*}

\end{proof}
\begin{lem}\label{classical}
Let $<\alpha_n,\beta_n>$ be generators for Schottky groups $\G_n$ in the upper-half space model $\mathbb{H}^3$ 
with $|\tr(\beta_n)|\to\infty$ and 
$\alpha_n=\bigl(\begin{smallmatrix} \lambda_n & 0\\0 & \lambda^{-1}_n\end{smallmatrix}\bigr)$, $|\lambda_n|>1$. Suppose one of the following set of conditions
holds:\par
there exists $\Lambda>1$, such that for large $n,$ we have $|\lambda_n|<\Lambda$ and, 
\begin{itemize}
\item
$|\lambda_n|^{-1}<|z_{\beta_n,l}|\le|z_{\beta_n,u}|<|\lambda_n|,$ and
\item 
\[\liminf_n\left\{\frac{1}{(|z_{\beta_n,u}|-|\lambda_n|)|\tr(\beta_n)|},\frac{1}{(|z_{\beta_n,l}|-|\lambda_n|^{-1})|\tr(\beta_n)|}
\right\}=0,\]
\end{itemize}
or there exists $\kappa>0$ and for large $n,$ we have $|\lambda_n|>\kappa$ and,
\begin{itemize}
\item
$\kappa^{-1}<|z_{\beta_n,l}|\le|z_{\beta_n,u}|<\kappa,$ and
\item 
\[\liminf_n\left\{\frac{1}{(|z_{\beta_n,u}|-\kappa)|\tr(\beta_n)|},\frac{1}{(|z_{\beta_n,l}|-\kappa^{-1})|\tr(\beta_n)|}
\right\}=0,\]
\end{itemize}
then there exists a subsequence such that for $i$ large, $\pi^{-1}<\alpha_{n_i},\beta_{n_i}>\pi$ 
are classical generators for $\G_{n_i}$ in the unit ball model $\mathbb{B}$.
\end{lem}
\begin{proof}
Let us suppose there exists a subsequence $<\alpha_{n_i},\beta_{n_i}>$ that satisfies the first set of conditions.
First assume that for large $i,$ $|z_{\beta_{{n_i},u}}-z_{\beta_{{n_i},l}}|>\delta>0.$  \par
Let $r_i,\rho_i$ denote the Euclidean radii of the isometric spheres of $\pi^{-1}\alpha_{n_i}\pi$ and 
$\pi^{-1}\beta_{n_i}\pi$ respectively. 
Note that $4\cosh(T_{\pi^{-1}\beta_{n_i}\pi})\ge|\tr^2(\pi^{-1}\beta_{n_i}\pi)|,$ which implies
there exists $c'>0$ such that $e^{T_{\pi^{-1}\beta_{n_i}\pi}}\ge c'|\tr^2(\pi^{-1}\beta_{n_i}\pi)|.$
Since $\rho^{-1}_i=\cosh\dis(o,\mathcal{L}_{\pi^{-1}\beta_{n_i}\pi})\sinh(\frac{1}{2}T_{\pi^{-1}\beta_{n_i}\pi})$ 
(\cite{Berdon} p$175$), 
we have $\rho^{-1}_i\ge\sinh(\frac{1}{2}T_{\pi^{-1}\beta_{n_i}\pi}),$ so for large $i$ there exists $c>0$ such that
$\rho^{-1}_i\ge ce^{\frac{1}{2}T_{\pi^{-1}\beta_{n_i}\pi}}.$
Hence 
there exists $\delta_1>0$ such that $\rho_i\le\delta_1|\tr(\beta_{n_i})|^{-1}$ for large $i.$
Since for $z,w\in\mathbb{C}$,
\[|\pi^{-1}(z)'|=\frac{2}{|z+{\bf e}|^2}\quad\quad\text{and},\]
\[|\pi^{-1}(z)-\pi^{-1}(w)|=|\pi^{-1}(z)'|^{1/2}|\pi^{-1}(w)'|^{1/2}|z-w|.\]
This implies for large $i,$ and $x_i\in C_{|z_{\beta_{n_i},u}|}$ and $y_i\in C_{|\lambda_{n_i}|},$
\begin{eqnarray*}
\frac{\rho_i}{|\pi^{-1}x_i-\pi^{-1}y_i|}&\le& \frac{2\delta_1|x_i+{\bf e}||y_i+{\bf e}|}{|\tr(\beta_{n_i})||x_i-y_i|}.
\end{eqnarray*}
Since $|\lambda_{n_i}|<\Lambda,$ there exists $\delta_2>0,$ such that $|x_i+{\bf e}||y_i+{\bf e}|<\delta_2$, and
\begin{eqnarray*}
\lim_i\frac{\rho_i}{|\pi^{-1}x_i-\pi^{-1}y_i|}
&\le&\lim_i\frac{2\delta_1\delta_2}{|\tr(\beta_{n_i})||x_i-y_i|}=0.
\end{eqnarray*}
Similarly there exists $\delta_3>0,$ such that for $w_i\in C_{|z_{\beta_{n_i},l}|}$ and $z_i\in C_{|\lambda_{n_i}|^{-1}},$
\begin{eqnarray*}
\lim_i\frac{\rho_i}{|\pi^{-1}w_i-\pi^{-1}z_i|}
&\le&\lim_i\frac{2\delta_1\delta_3}{|\tr(\beta_{n_i})||w_i-z_i|}=0.
\end{eqnarray*}
Hence it follows that for large $i,$ and Proposition \ref{proj}, the isometric spheres 
$S_{\pi^{-1}\alpha_{n_i}\pi,r_i},S_{\pi^{-1}\alpha^{-1}_{n_i}\pi,r_i},$
$S_{\pi^{-1}\beta_{n_i}\pi,\rho_i},S_{\pi^{-1}\beta^{-1}_{n_i}\pi,\rho_i}$ are disjoint.
\begin{rem}\label{classical-rem-2}
Note that if we don't assume that $|\lambda_{n_i}|<\Lambda$ then we don't have bounds on $|x_i+{\bf e}||y_i+{\bf e}|.$ 
However, since $|x_i+{\bf e}||y_i+{\bf e}|\le (|z_{\beta_{n_i},u}|+1)(|\lambda_{n_i}|+1)$ for $x_i\in C_{|z_{\beta_{n_i},u}|}, y_i\in C_{|\lambda_{n_i}|}$ and
$|x_i+{\bf e}||y_i+{\bf e}|\le (|z_{\beta_{n_i},l}|+1)(|\lambda_{n_i}|^{-1}+1)$ for $x_i\in C_{|z_{\beta_{n_i},l}|}, y_i\in C_{|\lambda_{n_i}|^{-1}}.$ 
Hence we can state the condition as follows,
\begin{itemize}
\item
$|\lambda_n|^{-1}<|z_{\beta_n,l}|\le|z_{\beta_n,u}|<|\lambda_n|,$ and
\item 
\[\liminf_n\left\{\frac{(|z_{\beta_n,u}|+1)(|\lambda_n|+1)}{(|z_{\beta_n,u}|-|\lambda_n|)|\tr(\beta_n)|},
\frac{(|z_{\beta_n,l}|+1)(|\lambda_n|^{-1}+1)}{(|z_{\beta_n,l}|-|\lambda_n|^{-1})|\tr(\beta_n)|}
\right\}=0,\]
\end{itemize}
\end{rem}
Next let us assume that $|z_{\beta_{n_i},u}-z_{\beta_{n_i},l}|\to 0.$ 
Under this assumption we can do a much stronger estimate of the lower bounds of $\cosh\dis(j,\mathcal{L}_{\beta_{n_i}}),$ distance between 
the point $j$ on the vertical $j$-axis and the axis of $\beta_{n_i}$ in the $\mathbb{H}^3.$ Note however that a weaker lower bounds is
sufficient in our case.\par
Recall that given any two points $h_1=(z_1,\theta_1),h_2=(z_2,\theta_2)\in\mathbb{H}^3$ the hyperbolic distance is given by,
\[\cosh\dis(h_1,h_2)=\frac{|z_1-z_2|^2+|\theta_1-\theta_2|^2}{2\theta_1\theta_2}+1.\]
By $|z_{\beta_{n_i},u}-z_{\beta_{n_i},l}|\to 0,$ also 
$\frac{1}{\Lambda}<|z_{\beta_{n_i},l}|\le|z_{\beta_{n_i},u}|<\Lambda$ we can estimate $\cosh\dis(j,\mathcal{L}_{\beta_{n_i}})$ 
by using the above formula for $(z_1,\theta_1)\in\mathcal{L}_{\alpha_{n_i}}$ and $(z_2,\theta_2)\in\mathcal{L}_{\beta_{n_i}}.$
Since for large $i$ we have, $|z_1-z_2|\ge|z_{\beta_{n_i},l}|,$
$|\theta_1-\theta_2|\ge||z_{\beta_{n_i},l}|-\frac{1}{2}(|z_{\beta_{n_i},u}-z_{\beta_{n_i},l}|)|,$ 
$2\theta_1\theta_2\le |z_{\beta_{n_i},u}||z_{\beta_{n_i},u}-z_{\beta_{n_i},l}|.$
Hence for large $i$ we have,
\[\cosh\dis(j,\mathcal{L}_{\beta_{n_i}})\ge\frac{|z_{\beta_{n_i},l}|^2+(|z_{\beta_{n_i},l}|-\frac{1}{2}(|z_{\beta_{n_i},u}-z_{\beta_{n_i},l}|))^2}
{|z_{\beta_{n_i},u}||z_{\beta_{n_i},u}-z_{\beta_{n_i},l}|}+1.\]
Since $|\Lambda|^{-1}<|z_{\beta_{n_i},l}|\le|z_{\beta_{n_i},u}|<\Lambda$ and $|z_{\beta_{n_i},u}-z_{\beta_{n_i},l}|\to 0,$  we have for large $i$ 
there exits $\sigma>0$
such that 
\[\cosh\dis(j,\mathcal{L}_{\beta_{n_i}})\ge\frac{\sigma}{|z_{\beta_{n_i},u}-z_{\beta_{n_i},l}|}.\]
Also by $\Lambda^{-1}<|z_{\beta_{n_i},l}|\le|z_{\beta_{n_i},u}|<\Lambda$ and $|\lambda_{n_i}|<|\Lambda|,$ there exists $\sigma'>0$ such that
$|\pi^{-1}(z_{\beta_{n_i},u})'||\pi^{-1}(z_{\beta_{n_i},l})'|>\sigma'.$
Since,
\[ |\pi^{-1}z_{\beta_{n_i},u}-\pi^{-1}z_{\beta_{n_i},l}|=|\pi^{-1}(z_{\beta_{n_i},u})'||\pi^{-1}(z_{\beta_{n_i},l})'|
|z_{\beta_{n_i},u}-z_{\beta_{n_i},l}|,\] we have
\[|\pi^{-1}z_{\beta_{n_i},u}-\pi^{-1}z_{\beta_{n_i},l}|\ge\sigma'|z_{\beta_{n_i},u}-z_{\beta_{n_i},l}|,\] 
From $\rho^{-1}_i=\cosh\dis(o,\mathcal{L}_{\pi^{-1}\beta_{n_i}\pi})\sinh(\frac{1}{2}T_{\pi^{-1}\beta_{n_i}\pi})$ and
the above estimates, implies that for $i$ large, there exists
$\delta_4>0$ such that $\rho_i\le\delta_4|\pi^{-1}z_{\beta_{n_i},u}-\pi^{-1}z_{\beta_{n_i},l}||\tr(\beta_{n_i})|^{-1}.$
Hence there exists $\delta_5>0$ such that for $x_i\in C_{|z_{\beta_{n_i},u}|}, y_i\in C_{|\lambda_{n_i}|},$
\begin{eqnarray*}
\lim_i\frac{\rho_i}{|\pi^{-1}x_i-\pi^{-1}y_i|}
&\le&\lim_i\frac{\delta_5|\pi^{-1}z_{\beta_{n_i},u}-\pi^{-1}z_{\beta_{n_i},l}|}{|\tr(\beta_{n_i})||x_i-y_i|}=0.
\end{eqnarray*}
Similarly there exists $\delta_6>0,$ such that for $w_i\in C_{|z_{\beta_{n_i},l}|}$ and $z_i\in C_{|\lambda_{n_i}|^{-1}},$
\begin{eqnarray*}
\lim_i\frac{\rho_i}{|\pi^{-1}w_i-\pi^{-1}z_i|}
&\le&\lim_i\frac{\delta_6|\pi^{-1}z_{\beta_{n_i},u}-\pi^{-1}z_{\beta_{n_i},l}|}{|\tr(\beta_{n_i})||w_i-z_i|}=0.
\end{eqnarray*}
From these estimates and Proposition \ref{proj} we have for sufficiently large $i,$ $S_{\pi^{-1}\alpha_{n_i}\pi,r_i},$ $S_{\pi^{-1}\alpha^{-1}_{n_i}\pi,r_i}$
are disjoint from $S_{\pi^{-1}\beta_{n_i}\pi,\rho_i},$ $S_{\pi^{-1}\beta^{-1}_{n_i}\pi,\rho_i},$ and since $|\tr(\beta_{n_i})|\to\infty$
implies $S_{\pi^{-1}\beta_{n_i}\pi,\rho_i}$ and $S_{\pi^{-1}\beta^{-1}_{n_i}\pi,\rho_i}$ are disjoint when $i$ is large, we have the first part 
of the lemma.\par
The second part of the lemma can be proved in the same way.
\end{proof}
\begin{rem}\label{classical-rem}
Note that in the course of the proof we see that if $|\lambda_n|\to 1$ then we can weaken the first set of conditions in the above lemma to:
\begin{itemize}
\item
$|\lambda_n|^{-1}<|z_{\beta_n,l}|\le|z_{\beta_n,u}|<|\lambda_n|,$ and
\item 
\[\liminf_n\left\{\frac{|z_{\beta_n,l}-z_{\beta_n,u}|}{(|z_{\beta_n,u}|-|\lambda_n|)|\tr(\beta_n)|},\frac{|z_{\beta_n,l}-z_{\beta_n,u}|}{(|z_{\beta_n,l}|-|\lambda_n|^{-1})|\tr(\beta_n)|}
\right\}=0,\]
\end{itemize}
\end{rem}
\section{Schottky Subspaces $\mf{J}_{k}(\tau)$}
This section is devoted to proving Theorem \ref{t-space} by utilizing results established in the previous sections.\par
For $\tau>0$, define $\mf{J}_k(\tau):=\{\G\in\mf{J}_{k}| Z_\G>\tau\}$. Recall that $\mf{J}_{k}$ denotes set of all Schottky groups
of rank $k$. 

\begin{thm}\label{t-space}
Let $\mf{J}_2$ be the set of all $2$-generated Schottky groups. 
For each $\tau>0$ there exists a $\nu>0$ such that 
$\{\G\in\mf{J}_2(\tau)| D_\G\le\nu\}\subset\mf{J}_{k,o}$.
\end{thm}
\begin{proof}
We prove by contradiction. Assume there exists a sequence $\{\G_n\}\subset\mf{J}_2(\tau)$
of nonclassical Schottky groups with $D_n\to 0$. By passing to subsequence, we may assume
$D_n\to 0$ monotonically.
Set $\G_n=<\alpha_n,\beta_n>$ with
$Z_{<\alpha_n,\beta_n>}>\tau$. We arrange the generators so that 
$|\tr(\alpha_n)|\le |\tr(\beta_n)|$. There are two possibilities: $(I)$ There exists a subsequence
such that $|\tr(\alpha_{n_i})|\to\infty$, $(II)$ $|\tr(\alpha_n)|<M$, for some $M>0.$\par
Case $(I)$ is trivial. Since both $|\tr(\alpha_{n_i})|,|\tr(\beta_{n_i})|\to \infty$ as $n\to\infty$,
it follows from $Z_{<\alpha_n,\beta_n>}>\tau$, there must exists $N$ such that
$<\alpha_n,\beta_n>$ becomes classical Schottky groups for $n>N$. A contradiction.\par
Now we consider case $(II)$. We work in upper space model $\mathbb{H}^3$. Conjugate $<\alpha_n,\beta_n>$ by 
a Mobius transformation into 
$\alpha_n=\bigl(\begin{smallmatrix} \lambda_n & 0\\0 & \lambda^{-1}_n\end{smallmatrix}\bigr)$
with $|\lambda_n|>1$. Denote 
$\beta_n=\bigl(\begin{smallmatrix} a_n & b_n\\c_n & d_n\end{smallmatrix}\bigr)$. Since $|\tr(\alpha_n)|<M$
implies $|\lambda_n|<M'$ for some $M'>0,$ it follows from Proposition \ref{trace} and $D_n\to 0$, we have $|\tr(\beta_n)|\to\infty.$
In addition, by conjugation with Mobius transformations, we can assume $\beta_n$ have $\eta_{\beta_n}=1.$
By replacing $\beta_n$ with $\beta^{-1}_n$ if necessary, we can assume $|\zeta_{\beta_n}|\le|\eta_{\beta_n}|.$\par
Since $Z_{<\alpha_n,\beta_n>}>\tau$, there exists $\Delta_1,\Delta_2,\Delta_3,\Delta_4>0$ such that, 
$\Delta_1<|z_{\beta_n,l}|\le|z_{\beta_n,u}|<\Delta_2,$ and
$\Delta_3<|z_{\beta_n,l}-z_{\beta_n,u}|<\Delta_4.$
It follows from Lemma \ref{fix-trace}, $\Delta_1<|\zeta_{\beta_n}|\le|\eta_{\beta_n}|<\Delta_2$, and
$\Delta_3<\lim_n|\zeta_{\beta_n}-\eta_{\beta_n}|<\Delta_4.$
\par
For each $n$, choose integers $k_n$ such that:
\[1\le|\zeta_{\beta_n}\lambda^{2k_n}_n|<|\lambda^2_n|.\]
We consider the generating set
$<\alpha_n,\alpha^{k_n}_n\beta_n>$. \par
By passing to subsequence if necessary, there are three cases that need to be considered 
(To simplify the notation, we denote subsequences by the same index notation): 
\begin{itemize}
\item[$(A)$]
$|\zeta_{\beta_n}\lambda^{2k_n}_n|\to 1,$ 
\item[$(B)$]
$||\zeta_{\beta_n}\lambda^{2k_n}_n|-|\lambda_n|^2|\to 0,$
\item[$(C)$]
case $(A),(B)$ do not occur. 
\end{itemize}
\subsection{Case $(C)$} 
We use the same notation index for subsequences. Since cases $(A),(B)$ do not occur, for large $n$ there exists 
$1<|\lambda|<|\lambda_n|$, $\sigma<1$, such that  
$|\zeta_n\lambda^{2k_n}_n|\to \sigma|\lambda|^2$.
Let $\psi_n$ be the Mobius transformation that fixes $\{0,\infty\}$ 
defined by $\psi(x)=\frac{x}{\sqrt{\sigma}|\lambda|}$, $x\in\mathbb{C}.$
Then $|\eta_{\psi\alpha^{k_n}_n\beta_n\psi^{-1}}|\to\frac{1}{\sqrt{\sigma}|\lambda|}$
and $|\zeta_{\psi\alpha^{k_n}_n\beta_n\psi^{-1}}|\to\sqrt{\sigma}|\lambda|.$
It follows from Lemma \ref{fix-trace}, for $n$ large, 
\[\max\{|z_{\psi\alpha^{k_n}_n\beta_n\psi^{-1},\mp}-\eta_{\psi\alpha^{k_n}_n\beta_n\psi^{-1}}|,
|z_{\psi\alpha^{k_n}_n\beta_n\psi^{-1},\pm}-\zeta_{\psi\alpha^{k_n}_n\beta_n\psi^{-1}}|\}
<\frac{\rho}{|\tr(\beta_n)|}.\]
Hence there exists $\rho', \rho'',\rho'''>0$ such that,
\[\left|z_{\psi\alpha^{k_n}_n\beta_n\psi^{-1},+}-z_{\psi\alpha^{k_n}_n\beta_n\psi^{-1},-}\right|>
\left||\sqrt{\sigma}|\lambda|-\frac{1}{\sqrt{\sigma}}|\lambda|\right|-\frac{\rho'}{|\tr(\beta_n)|},\]
and
\[\left|\frac{1}{z_{\psi\alpha^{k_n}_n\beta_n\psi^{-1},+}}-\frac{1}{z_{\psi\alpha^{k_n}_n\beta_n\psi^{-1},-}}\right|>
\rho''\left||\sqrt{\sigma}|\lambda|-\frac{1}{\sqrt{\sigma}}|\lambda|\right|-\frac{\rho'''}{|\tr(\beta_n)|}.\]
This implies that there exist $\Delta>0$ such that $Z_{\psi\alpha^{k_n}_n\beta_n\psi^{-1}}>\Delta.$
Hence by applying Proposition \ref{trace} to the generators $<\psi\alpha_n\psi^{-1},\psi\alpha^{k_n}_n\beta_n\psi^{-1}>$,
implies $|\tr(\psi\alpha^{k_n}_n\beta_n\psi^{-1})|\to\infty.$\par
Set $\kappa=\sqrt{\frac{\sigma+1}{2}}|\lambda|.$ Then for sufficiently large $n$ we have 
$\kappa<|\lambda_n|,$ $\kappa^{-1}<|z_{\psi\alpha^{k_n}_n\beta_n\psi^{-1},l}|\le
|z_{\psi\alpha^{k_n}_n\beta_n\psi^{-1},u}|<\kappa.$
And obviously,
\[ \lim_n\frac{1}{(\kappa-|z_{\psi\alpha^{k_n}_n\beta_n\psi^{-1},u}|)|\tr(\psi\alpha^{k_n}_n\beta_n\psi^{-1})|}=0,\]
\[\lim_n\frac{1}{(|z_{\psi\alpha^{k_n}_n\beta_n\psi^{-1},l}|-\kappa^{-1})|\tr(\psi\alpha^{k_n}_n\beta_n\psi^{-1})|}=0.\]
Therefore, $<\psi\alpha_n\psi^{-1},\psi\alpha^{k_n}_n\beta_n\psi^{-1}>$ satisfies the second set of conditions of
Lemma \ref{classical}, and so by Lemma \ref{classical}, these will be classical generators for large $n$, a contradiction.

\subsection{case $(A)$}
By passing to a subsequence if necessary,
we have two possibilities: 
\begin{itemize}
\item[$(A_1)$] 
$|\lambda_n|^2-1$ is monotonically decreasing to $0$. 
\item[$(A_2)$]
There exists $\lambda>1$, such that
$|\lambda_n|\ge|\lambda|$ for large $n$.
\end{itemize}
\subsubsection{$(A_1)$} 
Here we either have 
\[(i)\quad\quad\limsup_n\left||\zeta_{\alpha^{k_n}_n\beta_n}|-|\lambda_n|^2\right|
\left|\zeta_{\alpha^{k_n}_n\beta_n}-1\right|^{-1}<\infty \]
or
\[(ii)\quad\quad\liminf_n\left||\zeta_{\alpha^{k_n}_n\beta_n}|-|\lambda_n|^2\right|\left|\zeta_{\alpha^{k_n}_n\beta_n}-1\right|^{-1}\to\infty.\] 
If $\alpha^{k_n}_n\beta_n$ satisfies $(ii)$, we
conjugate $\alpha^{k_n}_n\beta_n$ with Mobius transformation $\psi_n$ defined by, 
$\psi_n(x)=\frac{x}{\zeta_{\alpha^{k_n}_n\beta_n}}.$
Consider $(\psi_n\alpha^{k_n}_n\beta_n\alpha^{l_n}_n\psi^{-1}_n)^{-1},$ take $l_n=-1$. 
Since by factor out $\zeta_{\alpha^{k_n}_n\beta_n}^2\lambda_n^{-2}$ in $(ii)$,
\[\liminf_n |\zeta_{\alpha^{k_n}_n\beta_n}|^2|\lambda_n|^{-2}\left||1-|\zeta_{\alpha^{k_n}_n\beta_n}|^{-1}|\lambda_n|^2\right|\left|\zeta^{-1}_{\alpha^{k_n}_n\beta_n}\lambda^2_n-\lambda_n^2\right|^{-1}\to\infty.\] 
\[\liminf_n \left||1-|\zeta_{\alpha^{k_n}_n\beta_n}^{-1}\lambda^2_n|\right|\left|\zeta^{-1}_{\alpha^{k_n}_n\beta_n}\lambda^2_n-\lambda_n^2\right|^{-1}\to\infty.\]
Since $\zeta_{\alpha^{k_n}_n\beta_n}^{-1}\lambda^2_n= \eta_{\alpha^{k_n}_n\beta_n\alpha^{-1}_n}$ and 
$\zeta_{(\psi_n\alpha^{k_n}_n\beta_n\alpha_n^{-1}\psi^{-1}_n)^{-1}} =\eta_{\psi_n\alpha^{k_n}_n\beta_n\alpha^{-1}_n\psi^{-1}_n}$ we have,
\[\liminf_n \left||1-|\eta_{\alpha^{k_n}_n\beta_n\alpha^{-1}_n}|^2\right|\left|\eta_{\alpha^{k_n}_n\beta_n\alpha_n^{-1}}-\lambda_n^2\right|^{-1}\to\infty,\] 
giving,
\[\limsup_n\left||\zeta_{(\psi_n\alpha^{k_n}_n\beta_n\alpha_n^{-1}\psi^{-1}_n)^{-1}}|-|\lambda_n|^2\right|
\left|\zeta_{(\psi_n\alpha^{k_n}_n\beta_n\alpha_n^{-1}\psi^{-1}_n)^{-1}}-1\right|^{-1}<\infty.\]
The generator $(\psi_n\alpha^{k_n}_n\beta_n\alpha^{-1}_n\psi^{-1}_n)^{-1}$ satisfies $(i)$. 
Hence replacing the generators if necessary
we can always assume the generators satisfy $(i)$. And without lost of generality we will assume that
$<\alpha_n,\alpha^{k_n}_n\beta_n>$ satisfies $(i).$\par
Consider $(i).$\par
In this case, we have either:
\[(i_1)\quad\quad\limsup_n\left||\zeta_{\alpha^{k_n}_n\beta_n}|-|\lambda_n|^2\right|
\left|\zeta_{\alpha^{k_n}_n\beta_n}-1\right|^{-1}>\delta>0, \]
or
\[(i_2)\quad\quad\limsup_n\left||\zeta_{\alpha^{k_n}_n\beta_n}|-|\lambda_n|^2\right|
\left|\zeta_{\alpha^{k_n}_n\beta_n}-1\right|^{-1}=0.\]
\paragraph{6.2.1.1} Consider $(i_1).$\par
\begin{lem}\label{axes}
There exists $c>0$ such that, 
\[\dis(\mathcal{L}_{\alpha_n},\mathcal{L}_{\alpha^{k_n}_n\beta_n})<\log(\frac{c}{|\lambda_n|^2-1}). \] 
\end{lem}
\begin{proof}
We first show that
$\frac{1}{|\tr(\beta_n)|(|\lambda_n|^2-1)}\to 0$. From Proposition \ref{trace} we have,
\begin{eqnarray*}
\lim\frac{1}{|\tr(\beta_n)|(|\lambda_n|^2-1)}&\le&\lim\rho\left(\frac{|\lambda_n|^{2D_n}-1}
{(|\lambda_n|^{2D_n}+3)(|\lambda_n|^2-1)^{2D_n}}\right)^{\frac{1}{2D_n}}
\end{eqnarray*}
and for large $n,$ we have $|\lambda_n|^{2D_n}-1<|\lambda_n|^2-1$ which implies that for some
$\rho'>0$,
\[\lim\frac{1}{|\tr(\beta_n)|(|\lambda_n|^2-1)}\le\lim\rho'(|\lambda_n|^2-1)^{\frac{1-2D_n}{2D_n}}=0.\]
It follows from Lemma \ref{fix-trace} and $\eta_{\alpha^{k_n}_n\beta_n}=1$, 
\[\left|\zeta_{\alpha^{k_n}_n\beta_n}-1\right|-
\rho'|\tr(\beta_n)|^{-1}\le\left|z_{\alpha^{k_n}_n\beta_n,-}-z_{\alpha^{k_n}_n\beta_n,+}\right|\le
 \left|\zeta_{\alpha^{k_n}_n\beta_n}-1\right|+
\rho''|\tr(\beta_n)|^{-1}\]
Since $1\le|\zeta_{\alpha^{k_n}_n\beta_n}|<|\lambda_n|^2,$ we have 
\[\frac{\left|z_{\alpha^{k_n}_n\beta_n,-}-z_{\alpha^{k_n}_n\beta_n,+}\right|}{|\lambda_n|^2-1}\ge 
\frac{\left|\zeta_{\alpha^{k_n}_n\beta_n}-1\right|}{|\lambda_n|^2-1}-\frac{\rho'}{|\tr(\beta_n)|(|\lambda_n|^2-1)}. \]
By the condition of $(i_{1})$ we have,
\[\frac{\left| \zeta_{\alpha^{k_n}_n\beta_n}-1\right|}{|\lambda_n|^2-1}=
\frac{|\zeta_{\alpha^{k_n}_n\beta_n}-1|}{|\lambda_n|^2-|\zeta_{\alpha^{k_n}_n\beta_n}|+
|\zeta_{\alpha^{k_n}_n\beta_n}|-1}>\frac{1}{M+1},\quad\quad\text{for some } M>0.\]
Hence for large $n$ there exists $\kappa>0$ such that,
\[\frac{\left|z_{\alpha^{k_n}_n\beta_n,-}-z_{\alpha^{k_n}_n\beta_n,+}\right|}{|\lambda_n|^2-1}>\frac{1}{M+1}-
\frac{\rho'}{|\tr(\beta_n)|(|\lambda_n|^2-1)}>\kappa.\] 
For the upper bounds we have, 
$|z_{\alpha^{k_n}_n\beta_n,-}-z_{\alpha^{k_n}_n\beta_n,+}|<2|\lambda_n|+\rho''|\tr(\beta_n)|^{-1}.$ 
Note that $\dis(\mathcal{L}_{\alpha_n},\mathcal{L}_{\alpha^{k_n}_n\beta_n})=
\inf\{\dis(h_1,h_2)|h_1\in\mathcal{L}_{\alpha_n},h_2\in\mathcal{L}_{\alpha^{k_n}_n\beta_n}\}.$ Set $h_j=(z_j,\theta_j); j=1,2$ then 
for a upper bound we can take 
\[(z_1,\theta_1)=(0, |z_{\alpha^{k_n}_n\beta_n,l}|+\frac{1}{2}|z_{\alpha^{k_n}_n\beta_n,u}-z_{\alpha^{k_n}_n\beta_n,l}|),\]
\[(z_2,\theta_2)=(\frac{1}{2}(z_{\alpha^{k_n}_n\beta_n,u}+z_{\alpha^{k_n}_n\beta_n,l}),
\frac{1}{2}|z_{\alpha^{k_n}_n\beta_n,u}-z_{\alpha^{k_n}_n\beta_n,l}|).\]
By Lemma \ref{fix-trace}, 
$1-\sigma_1|\tr(\beta_n)|^{-1}<|z_{\alpha^{k_n}_n\beta_n,l}|\le|z_{\alpha^{k_n}_n\beta_n,u}|\le|\lambda_n|^2+\sigma_2|\tr(\beta_n)|^{-1}$
and above estimates for $|z_{\alpha^{k_n}_n\beta_n,u}-z_{\alpha^{k_n}_n\beta_n,l}|$ we have,
\begin{eqnarray*}
\cosh\dis(\mathcal{L}_{\alpha_n},\mathcal{L}_{\alpha^{k_n}_n\beta_n})&\le&
\frac{\frac{1}{4}|z_{\alpha^{k_n}_n\beta_n,u}+z_{\alpha^{k_n}_n\beta_n,l}|^2+|z_{\alpha^{k_n}_n\beta_n,l}|^2}
{|z_{\alpha^{k_n}_n\beta_n,u}-z_{\alpha^{k_n}_n\beta_n,l}|(|z_{\alpha^{k_n}_n\beta_n,l}|+\frac{1}{2}|z_{\alpha^{k_n}_n\beta_n,u}-z_{\alpha^{k_n}_n\beta_n,l}|)}+1\\
&<&\frac{|\lambda_n|^2+\rho'''|\tr(\beta_n)|^{-1}
+\kappa(|\lambda_n|^2-1)+1}{\kappa(|\lambda_n|^2-1)},\quad \rho'''>0.
\end{eqnarray*}
This last inequality implies the Lemma.
\end{proof}
\begin{lem}\label{iso-ratio}
\[\lim_n\frac{1}{|\tr(\alpha^{k_n}_n\beta_n)|(|\lambda_n|^2-1)}=0\]
\end{lem}
\begin{proof}
It follows from Proposition \ref{trace} and  
Lemma \ref{axes}, there exists $\rho>0$ such that,
\[\left|\tr(\alpha^{k_n}_n\beta_n)\right|\ge\rho\left(\frac{(|\lambda_n|^{2D_n}+3)(|\lambda_n|^2-1)^{2D_n}}
{|\lambda_n|^{2D_n}-1}\right)^{\frac{1}{2D_n}}.\]
hence we have,
\[\lim\frac{1}{|\tr(\alpha^{k_n}_n\beta_n)|(|\lambda_n|^2-1)}\le
\lim\rho'\left(\frac{|\lambda_n|^{2D_n}-1}{(|\lambda_n|^{2D_n}+3)(|\lambda_n|^2-1)^{4D_n}}\right)^{\frac{1}{2D_n}}.\]
Since $|\lambda_n|^{2D_n}-1<|\lambda_n|^2-1$ for large $n$ we have,
\[\lim\frac{1}{|\tr(\alpha^{k_n}_n\beta_n)|(|\lambda_n|^2-1)}\le\lim\rho''(|\lambda_n|^2-1)^\frac{1-4D_n}{2D_n}=0.\]
\end{proof}
For large $n$ by condition $(i_1),$ we have $\delta(|\lambda_n|^2-1)<|\lambda_n|^2-|\zeta_{\alpha^{k_n}_n\beta_n}|$ then by Lemma \ref{fix-trace},
$||z_{\alpha^{k_n}_n\beta_n,u}|-|\zeta_{\alpha^{k_n}_n\beta_n}||\le
|z_{\alpha^{k_n}_n\beta_n,u}-\zeta_{\alpha^{k_n}_n\beta_n}|<\frac{\chi}{|\tr(\beta_n)|}$ for some $\chi>0$ we have,
\[|\lambda_n|^2-|z_{\alpha^{k_n}_n\beta_n,u}|>|\lambda_n|^2-|\zeta_{\alpha^{k_n}_n\beta_n}|-\frac{\chi}{|\tr(\beta_n)|}>
\delta(|\lambda_n|^2-1)-\frac{\chi}{|\tr(\beta_n)|}.\]
Set $\epsilon_n=(\delta(|\lambda_n|^2-1)-\frac{\chi}{|\tr(\beta_n)|}).$
Define Mobius transformations by $\psi_n(x)=(1+\frac{\epsilon_n}{2|z_{\alpha^{k_n}_n\beta_n,u}|})\lambda^{-1}_n(x).$ Then,
\begin{eqnarray*}
|\lambda_n|-|z_{\psi_n\alpha_n^{k_n}\beta_n\psi^{-1}_n,u}|&=&|\lambda_n|-(1+\frac{\epsilon_n}{2|z_{\alpha^{k_n}_n\beta_n,u}|})|\lambda_n|^{-1}
|z_{\alpha^{k_n}_n\beta_n,u}|\\
&=&(|\lambda_n|^2-|z_{\alpha^{k_n}_n\beta_n,u}|-\frac{\epsilon_n}{2}|)|\lambda_n|^{-1}\\
&>&(\epsilon_n-\frac{\epsilon_n}{2})|\lambda_n|^{-1}=\frac{\epsilon_n}{2}|\lambda_n|^{-1}.
\end{eqnarray*}
Also by Lemma \ref{fix-trace} and $\eta_{\alpha^{k_n}_n\beta_n}=1$ we have, $||z_{\alpha^{k_n}_n\beta_n,l}|-1|<\frac{\chi'}{|\tr(\beta_n)|}.$ 
This gives,
\begin{eqnarray*}
|z_{\psi_n\alpha_n^{k_n}\beta_n\psi^{-1}_n,l}|-|\lambda_n|^{-1}&=&(1+\frac{\epsilon_n}{2|z_{\alpha^{k_n}_n\beta_n,u}|})
|z_{\alpha_n^{k_n}\beta_n,l}||\lambda_n|^{-1}-|\lambda_n|^{-1}\\
&=&|\lambda_n|^{-1}\left(|z_{\alpha_n^{k_n}\beta_n,l}|-1+\frac{\epsilon_n|z_{\alpha_n^{k_n}\beta_n,l}|}{2|z_{\alpha^{k_n}_n\beta_n,u}|}\right)\\
&>&|\lambda_n|^{-1}\left(\frac{\epsilon_n|z_{\alpha_n^{k_n}\beta_n,l}|}{2|z_{\alpha^{k_n}_n\beta_n,u}|}-\frac{\chi'}{|\tr(\beta_n)|}\right)\\
&>&\frac{1}{|\lambda_n||\tr(\beta_n)|}\left(\frac{|\tr(\beta_n)|\epsilon_n|z_{\alpha_n^{k_n}\beta_n,l}|}{2|z_{\alpha^{k_n}_n\beta_n,u}|}-\chi'\right) 
\end{eqnarray*}
\[
>\frac{1}{|\lambda_n||\tr(\beta_n)|}\left(\frac{\delta|\tr(\beta_n)|(|\lambda_n|^2-1)|z_{\alpha_n^{k_n}\beta_n,l}|-\chi|z_{\alpha^{k_n}_n\beta_n,l}|}{2|z_{\alpha^{k_n}_n\beta_n,u}|}-\chi'\right)
\]
By Lemma \ref{iso-ratio} and above inequality we have, $|z_{\psi_n\alpha_n^{k_n}\beta_n\psi^{-1}_n,l}|-|\lambda_n|^{-1}>0$ and
$|\tr(\alpha^{k_n}_n\beta_n)|(|z_{\psi_n\alpha_n^{k_n}\beta_n\psi^{-1}_n,l}|-|\lambda_n|^{-1})\to\infty$ and 
$|\tr(\alpha^{k_n}_n\beta_n)|(|\lambda_n|-|z_{\psi_n\alpha_n^{k_n}\beta_n\psi^{-1}_n,u}|)\to\infty.$ Hence,
\begin{itemize}
\item
$|\lambda_n|^{-1}<|z_{\psi_n\alpha^{k_n}_n\beta_n\psi^{-1}_n,l}|\le|z_{\psi_n\alpha^{k_n}_n\beta_n\psi^{-1}_n,u}|<|\lambda_n|,$ and
\item 
\[\lim_n\left\{\frac{1}{(|z_{\psi_n\alpha^{k_n}_n\beta_n\psi_n^{-1},u}|-|\lambda_n|)|\tr(\alpha^{k_n}_n\beta_n)|},
\frac{1}{(|z_{\alpha^{k_n}_n\beta_n,l}|-|\lambda_n|^{-1})|\tr(\alpha^{k_n}_n\beta_n)|}
\right\}=0,\]
\end{itemize}
The generators 
$<\psi_n\alpha_n\psi^{-1}_n,\psi_n\alpha^{k_n}_n\beta_n\psi^{-1}_n>$ satisfies conditions of Lemma \ref{classical} for large $n.$
\paragraph{6.2.1.2} Consider $(i_2)$.\par
There exists $|\rho_{\alpha^{k_n-1}_n\beta_n}|\to 1$ such that $\zeta_{\alpha^{k_n-1}_n\beta_n}\rho_{\alpha^{k_n-1}\beta_n}=\eta_{\alpha^{k_n-1}_n\beta_n}=1$.
If $\limsup_n|\rho_{\alpha^{k_n-1}_n\beta_n}-1|>0$, then (with the same index notation for subsequence there exists a subsequence)
such that 
\[\liminf_n|\zeta_{\alpha^{k_n-1}_n\beta_n}-\eta_{\alpha^{k_n-1}_n\beta_n}|>0.\]
This implies by Lemma \ref{fix-trace}, $\liminf_n|z_{\alpha^{k_n-1}_n\beta_n,+}-z_{\alpha^{k_n-1}_n\beta_n,-}|>0.$ 
Hence by Proposition \ref{trace}, there exists $\rho>0$, such that for large $n$,
\[|\tr(\alpha^{k_n-1}_n\beta_n)|\ge\rho\left(\frac{|\lambda_n|^{2D_n}+3}{|\lambda_n|^{2D_n}-1}\right)^{\frac{1}{2D_n}}.\] In particular,
we have $|\tr(\alpha^{k_n-1}_n\beta_n)|\to\infty.$ Note 
\[\lim_n\frac{|\rho_{\alpha^{k_n-1}_n\beta_n}|-1}{|\lambda_n|^2-1}=0.\]
This can be seen as follows: since $|\lambda_n|^2-|\zeta_{\alpha^{k_n}_n\beta_n}|\le|\lambda_n|^2-1$ we have either, 
$\frac{|\lambda_n|^2-|\zeta_{\alpha^{k_n}_n\beta_n}|}{|\lambda_n|^2-1}\to 0$ or 
$1\ge\frac{|\lambda_n|^2-|\zeta_{\alpha^{k_n}_n\beta_n}|}{|\lambda_n|^2-1}>\epsilon>0.$\par
Assume that the latter inequality holds. This is equivalent to $(i_1)$ and we follows the same idea used in $(i_1).$
Set Mobius transformations $\psi_n(x)=\lambda_n^{-1}(1-\epsilon_n)^{-1}x,$ with $\epsilon_n=\frac{\epsilon(|\lambda_n|^2-1)}{2|\lambda_n|^2}.$ Then,
\begin{eqnarray*}
|\lambda_n|-|\zeta_{\psi_n\alpha^{k_n}_n\beta_n\psi^{-1}_n}|&=&|\lambda_n|-|\lambda_n^{-1}(1-\epsilon_n)^{-1}\zeta_{\alpha^{k_n}_n\beta_n}|\\
&=&|\lambda_n|^{-1}(1-\epsilon_n)^{-1}\left(|\lambda_n|^2(1-\epsilon_n)-|\zeta_{\alpha^{k_n}_n\beta_n}|\right)
\end{eqnarray*}
Since $\quad |\lambda_n|^2-|\zeta_{\alpha^{k_n}_n\beta_n}|>\epsilon(|\lambda_n|^2-1)$ we have,
\[|\lambda_n|-|\zeta_{\psi_n\alpha^{k_n}_n\beta_n\psi^{-1}_n}|> |\lambda_n|^{-1}(1-\epsilon_n)^{-1} (\epsilon(|\lambda_n|^2-1)-|\lambda_n|^2\epsilon_n)
=\frac{\epsilon(|\lambda_n|^2-1)}{2|\lambda_n|(1-\epsilon_n)}.\]
Since that $\epsilon_n\to 0,$ it follows from last inquality that for large $n$ we have,
$\frac{|\lambda_n|-|\zeta_{\psi_n\alpha^{k_n}_n\beta_n\psi^{-1}_n}|}{|\lambda_n|^2-1}>\epsilon'>0.$
And $\eta_{\psi_n\alpha_n^{k_n}\beta_n\psi^{-1}_n}=\lambda_n^{-1}(1-\epsilon_n)^{-1}$ we have, 
\[\frac{|\eta_{\psi_n\alpha_n^{k_n}\beta_n\psi^{-1}_n}|-|\lambda_n|^{-1}}{|\lambda|^2-1}=\frac{\epsilon}{2|\lambda_n|^3(1-\epsilon_n)}>\epsilon''>0.\]
Since $|\tr(\alpha^{k_n}_n\beta_n)(|\lambda_n|^2-1)|\to\infty$, it follows that for large $n$,
\begin{itemize}
\item
$|\lambda_n|^{-1}<|\eta_{\psi_n\alpha^{k_n}_n\beta_n\psi^{-1}_n}|\le|\zeta_{\psi_n\alpha^{k_n}_n\beta_n\psi^{-1}_n}|<|\lambda_n| ,$ and
\item 
$\lim_n\left\{\frac{1}{|\tr(\alpha^{k_n}_n\beta_n)|(|\eta_{\psi_n\alpha^{k_n}\beta_n\psi_n^{-1}}|-|\lambda_n|^{-1})},
\frac{1}{|\tr(\alpha^{k_n}_n\beta_n)|(|\zeta_{\psi_n\alpha^{k_n}_n\beta_n\psi_n^{-1}}|-|\lambda_n|)}\right\}=0.$
\end{itemize}
Hence by Lemma \ref{fix-trace}, $<\psi_n\alpha_n\psi_n^{-1},\psi_n\alpha^{k_n}_n\beta_n\psi_n^{-1}>$ satisfies Lemma \ref{classical}.\par
If the former holds then 
\begin{eqnarray*}
|\lambda_n|^2(1-|\zeta_{\alpha^{k_n}_n\beta_n}\lambda_n^{-2}|)&=&|\lambda_n|^2(1-|\zeta_{\alpha^{k_n-1}_n\beta_n}|)\\
&=&|\lambda_n|^2(|\zeta_{\alpha^{k_n-1}_n\beta_n}\rho_{\alpha^{k_n-1}\beta_n}| -|\zeta_{\alpha^{k_n-1}_n\beta_n}|)\\
&=&|\lambda_n|^2|\zeta_{\alpha^{k_n-1}_n\beta_n}|(|\rho_{\alpha^{k_n-1}\beta_n}|-1)
\end{eqnarray*}
The last equation implies that $\lim_n\frac{|\rho_{\alpha^{k_n-1}_n\beta_n}|-1}{|\lambda_n|^2-1}=0,$ 
and $|\tr(\alpha^{k_n-1}_n\beta_n)(|\lambda_n|^2-1)|\to\infty$,
it follows that for large $n$,
\begin{itemize}
\item
$|\lambda_n|^{-1}<|\zeta_{\alpha^{k_n-1}_n\beta_n}|\le|\eta_{\alpha^{k_n-1}_n\beta_n}|<|\lambda_n|,$ and
\item 
$\lim_n\left\{\frac{1}{|\tr(\alpha^{k_n-1}_n\beta_n)|(|\eta_{\alpha^{k_n-1}\beta_n}|-|\lambda_n|)},
\frac{1}{|\tr(\alpha^{k_n-1}_n\beta_n)|(|\zeta_{\alpha^{k_n-1}_n\beta_n}|-|\lambda_n|^{-1})}\right\}=0.$
\end{itemize}
Hence by Lemma \ref{fix-trace}, $<\alpha_n,\alpha^{k_n-1}_n\beta_n>$ satisfies Lemma \ref{classical}.\par
\noindent Consider the case that, $\rho_{\alpha^{k_n-1}_n\beta_n}\to 1$.\par
Let $\psi_n$ be the Mobius transformations given by,
\[\psi_n(x)=\frac{1}{\Delta_n}\frac{z_{\alpha^{k_n-1}_n\beta_n,-}x-z_{\alpha^{k_n-1}_n\beta_n,-}z_{\alpha^{k_n-1}_n\beta_n,+}}
{x-z_{\alpha^{k_n-1}_n\beta_n,-}}, \quad x\in\mathbb{C}.\]
where $\Delta_n=\frac{1}{z_{\alpha^{k_n-1}_n\beta_n,-}z_{\alpha^{k_n-1}_n\beta_n,+}-z^2_{\alpha^{k_n-1}_n\beta_n,-}}.$
Let $\phi_n$ be the Mobius transformations defined by,
$\phi_n(x)=\eta^{-1}_{\psi_n\alpha^{k_n}_n\beta_n\psi^{-1}_n}x.$ Then
$\zeta_{\phi_n\psi_n\alpha^{k_n}_n\beta_n\psi^{-1}_n\phi^{-1}_n}$ is given by,
\begin{eqnarray*}
\zeta_{\phi_n\psi_n\alpha^{k_n}_n\beta_n\psi^{-1}_n\phi^{-1}_n}=
\left(\frac{1-z_{\alpha^{k_n-1}_n\beta_n,-}}{1-z_{\alpha^{k_n-1}_n\beta_n,+}}\right)
\left(\frac{\alpha^{k_n}_n\beta_n(z_{\alpha^{k_n-1}_n\beta_n,-})-z_{\alpha^{k_n-1}_n\beta_n,+}}
{\alpha^{k_n}_n\beta_n(z_{\alpha^{k_n-1}_n\beta_n,+})-z_{\alpha^{k_n-1}_n\beta_n,-}}\right)
\end{eqnarray*}
To see this we do a simple computation. Set $\tilde{\alpha}_n=\phi_n\psi_n\alpha^{k_n-1}_n\beta_n\psi^{-1}_n\phi^{-1}_n$ and 
$\tilde{\beta_n}=\phi_n\psi_n\alpha^{k_n}_n\beta_n\psi^{-1}_n\phi^{-1}_n.$ 
Writing in matrix $\tilde{\beta}_n=\bigl(\begin{smallmatrix} \tilde{a}_n & \tilde{b}_n\\ \tilde{c_n} & \tilde{d}_n\end{smallmatrix}\bigr).$
Note that by our choice of $\phi_n$, we have $\eta_{\tilde{\beta}_n}=1,$ so $\tilde{c}_n=-\tilde{d}_n$ and $\zeta_{\tilde{\beta}_n}=\frac{\tilde{a}_n}{-\tilde{d}_n}.$
By straight forward matrix multiplications we have,
\begin{eqnarray*}
\tilde{a}_n=
-z^2_{\alpha^{k_n-1}_n\beta_n,-}\lambda_n^{k_n}a_n&-&z_{\alpha^{k_n-1}_n\beta_n,-}\lambda_n^{k_n}b_n+\\
&&z_{\alpha^{k_n-1}_n\beta_n,-}z_{\alpha^{k_n-1}_n\beta_n,+}
(z_{\alpha^{k_n-1}_n\beta_n,-}\lambda^{-k_n}_nc_n+\lambda_n^{-k_n}d_n)
\end{eqnarray*}
\begin{eqnarray*}
\tilde{a}_n&=&\frac{
\left(\frac{-z_{\alpha^{k_n-1}_n\beta_n,-}(\lambda_n^{k_n}a_nz_{\alpha^{k_n-1}_n\beta_n,-}+\lambda^{k_n}_nb_n)}
{z_{\alpha^{k_n-1}_n\beta_n,-}\lambda^{-k_n}_nc_n+\lambda^{-k_n}_nd_n}+z_{\alpha^{k_n-1}_n\beta_n,-}z_{\alpha^{k_n-1}_n\beta_n,+}\right)}
{(z_{\alpha^{k_n-1}_n\beta_n,-}\lambda^{-k_n}_nc_n+\lambda^{-k_n}_nd_n)^{-1}}\\
&=&\frac{-z_{\alpha^{k_n-1}_n\beta_n,-}(\alpha^{k_n}_n\beta_n(z_{\alpha^{k_n-1}_n\beta_n,-})-z_{\alpha^{k_n-1}_n\beta_n,+})}
{(z_{\alpha^{k_n-1}_n\beta_n,-}\lambda_n^{-k_n}c_n+\lambda^{-k_n}_nd_n)^{-1}}
\end{eqnarray*}
\begin{eqnarray*}
\tilde{d}_n=z_{\alpha^{k_n-1}_n\beta_n,-}z_{\alpha^{k_n-1}_n\beta_n,+}\lambda^{k_n}_na_n&+&z_{\alpha^{k_n-1}_n\beta_n,-}\lambda^{k_n}_nb_n-\\
&&z^2_{\alpha^{k_n-1}_n\beta_n,-}(z_{\alpha^{k_n-1}_n\beta_n,+}\lambda_n^{k_n}c_n\lambda_n^{-k_n}d_n)
\end{eqnarray*}
\begin{eqnarray*}
\tilde{d}_n&=&\frac{
\left(\frac{z_{\alpha^{k_n-1}_n\beta_n,-}(\lambda_n^{k_n}a_nz_{\alpha^{k_n-1}_n\beta_n,+}+\lambda^{k_n}_nb_n)}
{z_{\alpha^{k_n-1}_n\beta_n,+}\lambda^{-k_n}_nc_n+\lambda^{-k_n}_nd_n}-z^2_{\alpha^{k_n-1}_n\beta_n,-}\right)}
{(z_{\alpha^{k_n-1}_n\beta_n,+}\lambda^{-k_n}_nc_n+\lambda^{-k_n}_nd_n)^{-1}}\\
&=&\frac{z_{\alpha^{k_n-1}_n\beta_n,-}(\alpha^{k_n}_n\beta_n(z_{\alpha^{k_n-1}_n\beta_n,+})-z_{\alpha^{k_n-1}_n\beta_n,-})}
{(z_{\alpha^{k_n-1}_n\beta_n,+}\lambda_n^{-k_n}c_n+\lambda^{-k_n}_nd_n)^{-1}}
\end{eqnarray*}
Now by $\eta_n=1$ we have,
\[\frac{ z_{\alpha^{k_n-1}_n\beta_n,+} \lambda_n^{-k_n}c_n+\lambda^{-k_n}_nd_n}{z_{\alpha^{k_n-1}_n\beta_n,-}\lambda_n^{-k_n}c_n+\lambda^{-k_n}_nd_n}=
\frac{1-z_{\alpha^{k_n-1}_n\beta_n,+}}{1-z_{\alpha^{k_n-1}_n\beta_n,-}}.\]
Since $\zeta_{\tilde{\beta}_n} =\frac{\tilde{a}_n}{-\tilde{d}_n},$ the formula for $\zeta_{\tilde{\beta}_n}$ follows from above equations for
$\tilde{a}_n$ and $\tilde{d}_n.$\par
Let $\tilde{\lambda}_n$, denote the multiplier of $\tilde{\alpha}_n.$ \par
First, we need to get a estimate of the growth of $|\tr(\tilde{\beta}_n)|$ in terms of $|\tr(\beta_n)|.$ Note that 
$|\tr(\tilde{\beta}_n)|=|\tr(\alpha^{k_n}_n\beta_n)|.$  
\begin{rema}\label{tilde-trace}
There exists $\sigma,\sigma'>0$ such that 
\[\sigma'|\tr(\beta_n)|>|\tr(\tilde{\beta}_n)|>\sigma|\tr(\beta_n)|(|\lambda_n|^2-1)\quad\quad\text{for $n$ large}.\]
\end{rema}
In fact we only need the lower bound for $|\tr(\tilde{\beta}_n)|$.
\begin{proof} 
Since, 
\[ \frac{|\lambda_n|^2-1}{|\zeta_{\alpha_n^{k_n}\beta_n}-1|}+ \frac{|\zeta_{\alpha_n^{k_n}\beta_n}|-|\lambda_n|^2}{|\zeta_{\alpha_n^{k_n}\beta_n}-1|}
=\frac{|\zeta_{\alpha_n^{k_n}\beta_n}|-1}{|\zeta_{\alpha_n^{k_n}\beta_n}-1|}\le1,\]
and by condition $(i_2),$ $\frac{|\zeta_{\alpha_n^{k_n}\beta_n}|-|\lambda_n|^2}{|\zeta_{\alpha_n^{k_n}\beta_n}-1|}\to 0$ we have,
\[ \frac{|\zeta_{\alpha_n^{k_n}\beta_n}-1|}{|\lambda_n|^2-1}>\epsilon>0,\quad\quad\text{for large $n$}.\]
Recall $\eta_{\alpha^{k_n}_n\beta_n}=1$, and by Lemma \ref{fix-bound-2}, more precisely by equation $(1)$ in the proof of Lemma \ref{fix-bound-2} and $|\tr(\beta_n)|\asymp |c_n|,$
\[\frac{|\tr(\alpha^{k_n}_n\beta_n)|}{|\tr(\beta_n)|}>\sigma' |\zeta_{\alpha^{k_n}_n\beta_n}-\eta_{\alpha^{k_n}_n\beta_n}|
>\sigma'\epsilon(|\lambda_n|^2-1), \quad\text{for large $n$}.\]
The upper bound is trivial.
\end{proof}

\begin{lem}\label{m=0}
Assume that there exists a subsequence such that,
\[\lim_j\frac{|\zeta_{\tilde{\beta}_{n_j}}-1|}{|\tilde{\lambda}_{n_j}|^2-1}=0.\]
Then $<\tilde{\alpha}_{n_j},\tilde{\beta}_{n_j}>$ are classical Schottky generators for large $j.$
\end{lem}
\begin{proof}
We will show that $<\tilde{\alpha}_{n_j},\tilde{\beta}_{n_j}>$ satisfies conditions of Lemma \ref{classical} with Remark \ref{classical-rem}.\par
Since by Remark 6.A, $\sigma'|\tr(\beta_n)|>|\tr(\tilde{\beta}_n)|>\sigma|\tr(\beta_n)|(|\lambda_n|^2-1).$ In particular $|\tr(\tilde{\beta}_n)|\to\infty$, we have
\[\lim_j\frac{\left|z_{\tilde{\beta}_{n_j},+}-z_{\tilde{\beta}_{n_j},-}\right|}{\left|\zeta_{\tilde{\beta}_{n_j}}-1\right|}=
\lim_j\left|\frac{\sqrt{|\tr(\tilde{\beta}_{n_j})|^2-4}}{\tr(\tilde{\beta}_{n_j})}\right|=1\]
Since $|\tr(\tilde{\beta}_{n_j})|\to\infty$ implies that the isometric circles of $\tilde{\beta}_{n_j}$ are disjoint for large $j.$ Since
$\zeta_{\tilde{\beta}_{n_j}},\eta_{\tilde{\beta}_{n_j}}$ are centers of these isometric circles and so by disjointness, the radius of these
isometric circles must be
$<\frac{|\zeta_{\tilde{\beta}_{n_j}}-\eta_{\tilde{\beta}_{n_j}}|}{2}$ for large $j.$ In addition, each isometric circle contains one of the fixed point
$z_{\tilde{\beta}_{n_j},l}, z_{\tilde{\beta}_{n_j},u}.$ By our convention $z_{\tilde{\beta}_{n_j},l},z_{\tilde{\beta}_{n_j},u}$ are contained within
the isometric circles with centers $\zeta_{\tilde{\beta}_{n_j}},\eta_{\tilde{\beta}_{n_j}}$ respectively.
Note that $\eta_{\tilde{\beta}_{n_j}}=1.$ Hence for large $j$ we have,
\[\frac{|z_{\tilde{\beta}_{n_j},l}-\zeta_{\tilde{\beta}_{n_j}}|}{|\zeta_{\tilde{\beta}_{n_j}}-1|},
\frac{|z_{\tilde{\beta}_{n_j},u}-1|}{|\zeta_{\tilde{\beta}_{n_j}}-1|}<\frac{1}{2}.\]
From these bounds we have,
\begin{align*}
\lim_j\frac{|z_{\tilde{\beta}_{n_j},l}-z_{\tilde{\beta}_{n_j},u}|}{(|z_{\tilde{\beta}_{n_j},u}|-|\tilde{\lambda}_{n_j}|)|\tr(\tilde{\beta}_{n_j})|}
&=\lim_j\frac{\left|\zeta_{\tilde{\beta}_{n_j}}-1\right|}{\left|(|z_{\tilde{\beta}_{n_j},u}|-1)+(1-|\tilde{\lambda}_{n_j}|)\right|\left|\tr(\tilde{\beta}_{n_j})\right|}\\
&\le\lim_j\frac{1}{\left(\left|\frac{|\tilde{\lambda}_{n_j}|-1}{\zeta_{\tilde{\beta}_{n_j}}-1}\right|-\frac{1}{2}\right)|\tr(\tilde{\beta}_{n_j})|}=0
\end{align*}
Similarly we have,
\begin{equation*}
\begin{split}
\lim_j\frac{|z_{\tilde{\beta}_{n_j},l}-z_{\tilde{\beta}_{n_j},u}|}{(|z_{\tilde{\beta}_{n_j},l}|-|\tilde{\lambda}_{n_j}|^{-1})|\tr(\tilde{\beta}_{n_j})|}
\quad\quad\quad\quad\quad\quad\quad\quad\quad\quad\quad\quad\quad\quad\quad\quad\quad\quad\quad\text{        }\\
=\lim_j\frac{\left|\zeta_{\tilde{\beta}_{n_j}}-1\right|}{\left|(|z_{\tilde{\beta}_{n_j},u}|-|\zeta_{\tilde{\beta}_{n_j}}|)+
(|\zeta_{\tilde{\beta}_{n_j}}|-1)+(1-|\tilde{\lambda}_{n_j}|^{-1})\right|\left|\tr(\tilde{\beta}_{n_j})\right|}\quad\quad\quad\quad\quad\text{        }\\
\le\lim_j\frac{1}{\left(\frac{1}{|\tilde{\lambda}_{n_j}|}\left|\frac{|\tilde{\lambda}_{n_j}|-1}{\zeta_{\tilde{\beta}_{n_j}}-1}\right|
+\left|\frac{|\zeta_{\tilde{\beta}_{n_j}}|-1}{\zeta_{\tilde{\beta}_{n_j}}-1}\right| -\frac{1}{2}\right)|\tr(\tilde{\beta}_{n_j})|}=0
\quad\quad\quad\quad\quad\quad\quad\quad\quad\quad\text{        }
\end{split}
\end{equation*}
Hence by Remark \ref{classical-rem}, $<\tilde{\alpha}_{n_j},\tilde{\beta}_{n_j}>$ are classical Schottky generators for sufficiently large $j.$
\end{proof}

\begin{lem}\label{axis-m=0}
Assume that there exists a $M>0$ such that,
\[M<\frac{|\zeta_{\tilde{\beta}_{n}}-1|}{|\tilde{\lambda}_{n}|^2-1}.\]
Then there exists $\mathcal{N},\sigma>0$ such that, $|\tilde{\lambda}_n|^{2D_n}-1>\sigma|\tr(\tilde{\beta}_n)|^{\frac{2D_n}{2D_n-1}},$
for $n>\mathcal{N}.$
\end{lem}
\begin{proof}
Use matrix representation we can write, 
$\tilde{\beta}_n=\bigl(\begin{smallmatrix} \tilde{a}_n & \tilde{b}_n\\ \tilde{c_n} & \tilde{d}_n\end{smallmatrix}\bigr).$
Note that $\eta_{\tilde{\beta}_n}=1.$ So $\tilde{c}_n=-\tilde{d}_n$ we have,
$\tilde{\beta}_n=\bigl(\begin{smallmatrix} \tilde{a}_n & \tilde{b}_n\\ -\tilde{d_n} & \tilde{d}_n\end{smallmatrix}\bigr).$
By $|\zeta_{\tilde{\beta}_n}|\le 1$ and our assumption of the Lemma 
$M<\frac{|\zeta_{\tilde{\beta}_{n}}-1|}{|\tilde{\lambda}_{n}|^2-1},$ 
we have, 
\[ M(|\tilde{\lambda}_n^2|-1)<\left|\zeta_{\tilde{\beta}_n}-1\right|<M'.\]
By $|\tr(\tilde{\beta}_n)|\to\infty$ and as in the proof of Lemma \ref{m=0} we have,
\[\lim_n\frac{\left|z_{\tilde{\beta}_{n},+}-z_{\tilde{\beta}_{n},-}\right|}{\left|\zeta_{\tilde{\beta}_{n}}-1\right|}=
\lim_n\left|\frac{\sqrt{|\tr(\tilde{\beta}_{n})|^2-4}}{\tr(\tilde{\beta}_{n})}\right|=1\]
Since $|z_{\tilde{\beta}_n,+}-z_{\tilde{\beta}_n,-}|=\left|\tilde{d}_n\right|^{-1}\left|\sqrt{\tr^2{\tilde{\beta}_n}-4}\right|$
it follows that we have for some $\sigma, \sigma'>0$ and large $n$ that,
\[\sigma|\tr(\tilde{\beta}_n)|<|\tilde{d}_n|<\sigma'|\tr(\tilde{\beta}_n)|(|\tilde{\lambda}_n|^2-1)^{-1}.\]
Let $e$ be the Euler number. If there exists a subsequence such that $\lim_i|\zeta_{\tilde{\beta}_{n_i}}|=e^{-1},$ 
then $1-e^{-1}<|z_{\tilde{\beta}_{n_i},+}-z_{\tilde{\beta}_{n_i},-}|<1+e^{-1}$ for large $n.$ Otherwise we let $m_n>0$ be integers defined as,
\[ 1+\frac{1}{m_n+1}\le|\tilde{\lambda}_n|\le 1+\frac{1}{m_n}.\]
From this definition we have, $\lim_n|\tilde{\lambda}_n|^{m_n}=\lim_n(1+m_n^{-1})^{m_n}=e.$
Then there exists $N,\delta>0$ such that for $n>N$ we have,
\begin{eqnarray*}
|\zeta_{\tilde{\alpha}^{m_n}_n\tilde{\beta}_n}-\eta_{\tilde{\alpha}^{m_n}_n\tilde{\beta}_n}|&\le&
\left||\zeta_{\tilde{\beta}_n}||\tilde{\lambda}_n|^{m_n}+1\right|\\
&\le& \delta(e+1) \\
\text{and,}\\
|\zeta_{\tilde{\alpha}^{m_n}_n\tilde{\beta}_n}-\eta_{\tilde{\alpha}^{m_n}_n\tilde{\beta}_n}|&\ge&
\left||\zeta_{\tilde{\beta}_n}||\tilde{\lambda}_n|^{m_n}-1\right|\\
&\ge& \delta(e-1). 
\end{eqnarray*}
Hence it follows from Lemma \ref{fix-trace}, there exists $\kappa>0$ such that,
\[ \kappa^{-1}<|z_{\tilde{\alpha}^{m_n}_n\tilde{\beta}_n,+}-z_{\tilde{\alpha}^{m_n}_n\tilde{\beta}_n,-}|<\kappa.  \]  
Therefore by setting $m_n=0$ for the subsequence with $\lim_i|\zeta_{\tilde{\beta}_{n_i}}|=e^{-1},$ we can always assume that for large $n$
\[ \kappa^{-1}<|z_{\tilde{\alpha}^{m_n}_n\tilde{\beta}_n,+}-z_{\tilde{\alpha}^{m_n}_n\tilde{\beta}_n,-}|<\kappa.  \]  
Since
\[|z_{\tilde{\alpha}^{m_n}_n\tilde{\beta}_n,+}-z_{\tilde{\alpha}^{m_n}_n\tilde{\beta}_n,-}|=
\frac{\left|\sqrt{\tr^2(\tilde{\alpha}^{m_n}_n\tilde{\beta}_n)-4}\right|}{|\tilde{d}_n\tilde{\lambda}^{m_n}_n|},\]
we have
\[\sigma''|\tr(\tilde{\beta}_n)|<|\tr(\tilde{\alpha}^{m_n}_n\tilde{\beta}_n)|<\frac{\sigma'''|\tr(\tilde{\beta}_n)|}{(|\tilde{\lambda}_n|^2-1)}.\]
Since $\kappa^{-1}<|z_{\tilde{\alpha}^{m_n}_n\tilde{\beta}_n,+}-z_{\tilde{\alpha}^{m_n}_n\tilde{\beta}_n,-}|<\kappa,$ we have
$\dis(\mathcal{L}_{\tilde{\alpha}_n},\mathcal{L}_{\tilde{\alpha}^{m_n}_n\tilde{\beta}_n})<\delta$ for some $\delta>0.$
By Remark 4.1.A of Proposition \ref{trace} applied to $<\tilde{\alpha}_n,\tilde{\alpha}_n^{m_n}\tilde{\beta}_n>$ we have,
\[ |\tr(\tilde{\alpha}^{m_n}_n\tilde{\beta}_n)|>
\rho\left(\frac{|\tilde{\lambda}_n|^{2D_n}+3}{|\tilde{\lambda}_n|^{2D_n}-1}\right)^{\frac{1}{2D_n}}\]
and above bound for $|\tr(\tilde{\alpha}^{m_n}_n\tilde{\beta}_n)|$ we have,
\begin{align*}
|\tr(\tilde{\beta}_n)|&>\rho\sigma'''^{-1}(|\tilde{\lambda}_n|^2-1)\left(\frac{|\tilde{\lambda}_n|^{2D_n}+3}{|\tilde{\lambda}_n|^{2D_n}-1}\right)^{\frac{1}{2D_n}}>
\rho'\frac{|\tilde{\lambda}_n|^2-1}{\left(|\tilde{\lambda}_n|^{2D_n}-1\right)^{\frac{1}{2D_n}}}\\
&>\rho''\left(|\tilde{\lambda}_n|^{2D_n}-1\right)^{1-\frac{1}{2D_n}}
\end{align*}
The last inequality implies that,
\[|\tilde{\lambda}_n|^{2}-1>|\tilde{\lambda}_n|^{2D_n}-1>\rho'''|\tr(\tilde{\beta}_n)|^{\frac{2D_n}{2D_n-1}}.\]
\end{proof}
\begin{prop}\label{classical-i-2}
Assume $(i_2)$.
Suppose that there exists a $M>0$ such that,
\[M<\frac{|\zeta_{\tilde{\beta}_{n}}-1|}{|\tilde{\lambda}_{n}|^2-1}.\]
Then $<\alpha_n,\alpha^{k_n-1}\beta_n>$ are classical generators for large $n.$
\end{prop}
To prove Proposition \ref{classical-i-2} when $\limsup_n|\tr(\alpha^{k_n-1}_n\beta_n)|<\infty$ we use disjoint non-isometric circles 
for $\alpha^{k_n-1}_n\beta_n$ based on the following Lemma.
\begin{prop}\label{non-iso}
Given any loxodromic transformation $\gamma$ with fixed points $\not=0,\infty$ and
mutiplier $\lambda^2_\gamma,$
there exists disjoint circles $\mathcal{S}_{o,r},\mathcal{S}_{o',r'}$ of center $o$ radius $r$ and center $o'$ radius $r'$ respectively such that, 
\[\gamma(\text{interior}(\mathcal{S}_{o,r}))\cap\text{interior}(\mathcal{S}_{o',r'})=\emptyset,\quad\text{and}\quad
r+r'=|z_{\gamma,+}-z_{\gamma,-}|\frac{2|\lambda_\gamma|}{|\lambda_\gamma|^2-1}.\]
Note that since $|\lambda_\g|>1,$ so by this equality for $r+r'$ we have a upper bound as, $r+r'<|z_{\gamma,-}-z_{\gamma,+}|\frac{|\lambda_\g|+1}{|\lambda_\g|-1}.$
\end{prop}
\begin{proof}
We conjugate $\gamma$ into
Mobius transformation $\gamma'$ with fixed points $\{0,\infty\}.$  Consider circles $\mathcal{S}_{0,|\lambda_\gamma|^{-1}},\mathcal{S}_{0,|\lambda_\gamma|}.$
The Mobius transformation $\phi(x)=\frac{x-1}{x+1},$ maps the fixed points of $\gamma'$ which are $\{0,\infty\}$ to fixed points $\{-1,1\}$ respectively.
In addition it maps $\mathcal{S}_{0,|\lambda_\gamma|^{-1}},\mathcal{S}_{0,|\lambda_\gamma|}$ to $S_{z_1,r_1}, S_{z'_1,r'_1}$ respectively. Here we can use basic
formulas to determine $z_1,z'_1,r_1,r'_1$ (see page $91$ of \cite{MCD}). Explicitly by \cite{MCD} we have,
\begin{eqnarray*}
r_1&=&\left|\frac{-|\lambda_\g|^{-2}-1}{-|\lambda_\g|^{-2}+1}-\frac{|\lambda_\g|^{-1}-1}{|\lambda_\g|^{-1}+1}\right|\\
&=&\frac{2|\lambda_\g|}{|\lambda_\g|^2-1}\\
r'_1&=&\left|\frac{-|\lambda_\g|^2-1}{-|\lambda_\g|^2+1}-\frac{|\lambda_\g|-1}{|\lambda_\g|+1}\right|\\
&=&\frac{2|\lambda_\g|}{|\lambda_\g|^2-1}
\end{eqnarray*}
This gives, 
\[r_1+r'_1=\frac{4|\lambda_\g|}{|\lambda_\g|^2-1}.\]
The distance between the centers is,
\[|z_1-z'_1|=\left|\frac{-|\lambda_\g|^{-2}-1}{-|\lambda_\g|^{-2}+1}-\frac{-|\lambda_\g|^2-1}{-|\lambda_\g|^2+1}\right|
=2\frac{|\lambda_\g|^2+1}{|\lambda_\g|^2-1}.\]
Since $(|\lambda_\g|^2+1)-2|\lambda_\g|=(|\lambda_\g|-1)^2>0,$ implies $\mathcal{S}_{z_1,r_1},\mathcal{S}_{z'_1,r'_1}$ are disjoint.
By conjugating $\phi\gamma'\phi^{-1}$ with $\psi(x)=x+\frac{z_{\gamma,+}+z_{\gamma,-}}{z_{\gamma,+}-z_{\gamma,-}}$ we map the fixed points 
$\{-1,1\}$ to $\left\{\frac{2z_{\g,-}}{z_{\g,+}-z_{\g,-}},\frac{2z_{\g,+}}{z_{\g,+}-z_{\g,+}}\right\}.$ 
Because $\psi(x)$ is a translation (i.e euclidian isometry), the circles are mapped to $\mathcal{S}_{z_2,r_2},\mathcal{S}_{z'_2,r'_2}$
with same radius and preserves the disjointness. Finally conjugate $\psi\phi\g\phi^{-1}\psi^{-1}$ by $\theta(x)=\frac{z_{\g,+}-z_{\g,-}}{2}$
maps $\left\{\frac{2z_{\g,-}}{z_{\g,+}-z_{\g,-}},\frac{2z_{\g,+}}{z_{\g,+}-z_{\g,+}}\right\}$ to $\{z_{\g,-},z_{\g,+}\},$ and maps circles to
$\mathcal{S}_{z_3,r_3},\mathcal{S}_{z'_3,r'_3}.$ Note that $r_3,r'_3=|z_{\g,+}-z_{\g,-}|\frac{|\lambda_\g|}{|\lambda_\g|^2-1}$ and preserves the disjointness.
Note that $\g=\theta\psi\phi\g'\phi^{-1}\psi^{-1}\theta^{-1}.$
Hence we have the sum of the radius of the resulting disjointed circles given by,
\[r_3+r_3'=|z_{\gamma,-}-z_{\gamma,+}|\frac{2|\lambda_\g|}{|\lambda_\g|^2-1}\]
\end{proof}
\begin{proof}{(Proposition \ref{classical-i-2})}
First assume that $\limsup_n|\tr(\alpha^{k_n-1}_n\beta_n)|<\infty.$
Let $\mathcal{S}_{o_n,r_n},\mathcal{S}_{o'_n,r'_n}$ be the disjoint circles for $\alpha^{k_n-1}_n\beta_n$ given by Proposition \ref{non-iso}.  We will
show that $\lim_n\frac{r_n+r'_n}{|\lambda_n|^2-1}=0.$\par
Note that
\begin{align*}
\lim_n\frac{r_n+r'_n}{|\lambda_n|^2-1}&\le\lim_n|z_{\alpha^{k_n-1}_n\beta_n,+}-z_{\alpha^{k_n-1}_n\beta_n,-}|
\frac{2|\tilde{\lambda}_n|}{(|\tilde{\lambda}_n|^2-1)(|\lambda_n|^2-1)}\\
&\text{by $|\tr(\alpha^{k_n-1}_n\beta_n)|<C$ for some $C>0$ we have,}\\
&\le\lim_n|z_{\alpha^{k_n-1}_n\beta_n,+}-z_{\alpha^{k_n-1}_n\beta_n,-}|\frac{2C}{(|\tilde{\lambda}_n|^2-1)(|\lambda_n|^2-1)}
\end{align*}
Since $\limsup_n|\tr(\alpha^{k_n-1}\beta_n)|<\infty$ and $\limsup_n|\lambda_n|^{k_n-1}<\infty,$ we have by Lemma \ref{fix-bound-2},
\[|z_{\alpha^{k_n-1}_n\beta_n,+}-z_{\alpha^{k_n-1}_n\beta_n,-}|\asymp\frac{1}{|\tr(\beta_n)|}.\]
By Proposition \ref{trace}, 
\[|\tr(\beta_n)|^{2D_n}>\rho^{D_n}\frac{|\lambda_n|^{2D_n}+3}{|\lambda_n|^{2D_n}-1}>\frac{4\rho^{D_n}}{|\lambda_n|^2-1}.\]
The last inequality follows from $|\lambda_n|^{2D_n}-1<|\lambda_n|^2-1$ for large $n.$\par
By our assumptiom that $M<\frac{|\zeta_{\tilde{\beta}_{n}}-1|}{|\tilde{\lambda}_{n}|^2-1},$ and Lemma \ref{axis-m=0} we have,
\begin{align*}
\lim_n\frac{r_n+r'_n}{|\lambda_n|^2-1}&\le\lim_n (4\rho^{D_n})^{-1}M'|\tr(\beta_n)|^{-1}|\tr(\beta_n)|^{2D_n}|\tr(\tilde{\beta}_n)|^{\frac{2D_n}{1-2D_n}}\\
&\text{since $\quad |\tr(\tilde{\beta}_n)|<\sigma|\tr(\beta_n)|$ we have,}\\
&<\lim_n (4\rho^{D_n})^{-1}\sigma^{\frac{2D_n}{1-2D_n}}M'|\tr(\beta_n)|^{\frac{D_n(6-4D_n)-1}{1-2D_n}}=0
\end{align*}
Now the circles $\mathcal{S}_{o_n,r_n}$ contains one of $z_{\alpha^{k_n-1}_n\beta_n,-},z_{\alpha^{k_n-1}\beta_n,+}$ and $\mathcal{S}_{o'_n,r'_n}$ contains
the other fixed point, and since $\eta_{\alpha^{k_n-1}_n\beta_n}=1$ and $|\zeta_{\alpha^{k_n-1}_n\beta_n}|\to 1,$ 
it follows from Lemma \ref{fix-trace}, we must have $\mathcal{S}_{o_n,r_n}, \mathcal{S}_{o'_n,r'_n}$
contained in the region between $\frac{1}{|\lambda_n|}$ and $|\lambda_n|$ for large $n.$ Hence we have classical generators for large $n.$\par
Now if there exists a subsequence such that $|\tr(\alpha^{k_{n_i}-1}_{n_i}\beta_{n_i})|\to\infty,$ then $<\alpha_{n_i},\alpha^{k_{n_i}-1}\beta_{n_i}>$ satisfies
conditions of Lemma \ref{classical}.
\end{proof}

\subsubsection{$(A_2)$}
To prove $(A_2)$ we can follow the steps given in the proof of $(A_1),$ and do the approriate modifications. Some of the estimates will be 
simpler because $|\lambda_n|>\lambda$ and so estimates involving $(|\lambda_n|^2-1)^{-1}$ will hold trivially. However to avoid too much reproduction of the 
previous proof of $(A_1)$, we give here a alternative short cut proof of $(A_2)$ instead.\par
\begin{proof}{(of $(A_2)$)}
Suppose that there is a subsequence (use same index for subsequence) such that $|\tr(\alpha_n^{k_n}\beta_{n})|\to\infty.$ Then by our assumption of
$(A_2),$ $|\zeta_{\alpha^{k_n}_n\beta_n}|\to 1$ and $|\lambda_n|>\lambda>1$ we have for large $n$,  
\begin{itemize}
\item
$\lambda^{-1}< |\eta_{\alpha^{k_n}_n\beta_n}| \le|\zeta_{\alpha^{k_n}_n\beta_n}|<\lambda,$ and
\item 
$\lim_n\left\{\frac{1}{|\tr(\alpha^{k_n}_n\beta_n)|(|\eta_{\alpha^{k_n-1}\beta_n}|-\lambda^{-1})},
\frac{1}{|\tr(\alpha^{k_n}_n\beta_n)|(|\zeta_{\alpha^{k_n}_n\beta_n}|-\lambda_n)}\right\}=0.$
\end{itemize}
By Lemma \ref{fix-trace}, $<\alpha_n,\alpha^{k_n}_n\beta_n>$ satisfies the second set of conditions of Lemma \ref{classical}, hence
classical.\par
Otherwise we have $|\tr(\alpha^{k_n}_n\beta_n)|<C$ for some $C>0.$ Since $|\zeta_{\alpha^{k_n}_n\beta_n}|\to 1$ and $\eta_{\alpha^{k_n}_n\beta_n}=1,$
by Lemma \ref{fix-trace} we have $|z_{\alpha^{k_n}_n\beta_n,\pm}|\to 1.$ Now by Remark \ref{4A}.1.A and $|\tr(\alpha^{k_n}_n\beta_n)|<C$ 
we must have $\dis(\mathcal{L}_{\alpha_n},\mathcal{L}_{\alpha^{k_n}_n\beta_n})\to\infty.$ This implies that 
$|z_{\alpha^{k_n}_n\beta_n,+}-z_{\alpha^{k_n}_n\beta_n,-}|\to 0.$ More precisely we have,
\begin{lem}\label{delta-z}
Suppose $|z_{\alpha^{k_n}_n\beta_n,+}-z_{\alpha^{k_n}_n\beta_n,-}|\to 0$ and $|z_{\alpha^{k_n}_n\beta_n,\pm}|\to 1.$ 
Then there exists $\delta>0$ such that 
\[\dis(\mathcal{L}_{\alpha_n},\mathcal{L}_{\alpha^{k_n}_n\beta_n})<\log\left(\frac{\delta}{|z_{\alpha^{k_n}_n\beta_n,+}-z_{\alpha^{k_n}_n\beta_n,-}|}\right).\]
\end{lem}
\begin{proof}
By using hyperbolic distance formula, and since $|z_{\alpha^{k_n}_n\beta_n,+}-z_{\alpha^{k_n}_n\beta_n,-}|\to 0$ and
$|z_{\alpha^{k_n}_n\beta_n,\pm}|\to 1.$ we have for large $n$,
\begin{eqnarray*}
\cosh\dis(\mathcal{L}_{\alpha_n},\mathcal{L}_{\alpha^{k_n}_n\beta_n})&<&
\frac{|z_{\alpha^{k_n}_n\beta_n,u}|^2+\frac{1}{4}(|z_{\alpha^{k_n}_n\beta_n,u}+z_{\alpha^{k_n}_n\beta_n,l}|)^2}
{\frac{1}{2}(|z_{\alpha^{k_n}_n\beta_n,l}|+\frac{1}{2}|z_{\alpha^{k_n}_n\beta_n,u}-z_{\alpha^{k_n}_n\beta_n,l}|) 
(|z_{\alpha^{k_n}_n\beta_n,u}-z_{\alpha^{k_n}_n\beta_n,l}| )}+1\\
&<&\frac{\rho}{|z_{\alpha^{k_n}_n\beta_n,u}-z_{\alpha^{k_n}_n\beta_n,l}|}, \quad\text{for some $\rho>0$}
\end{eqnarray*}
\end{proof}
Let $\mathcal{S}_{o,r_n},\mathcal{S}_{o',r'_n}$ be the circles given by Proposition \ref{non-iso}.
\begin{prop}\label{delta-z-r}
$|\tr(\alpha^{k_n}_n\beta_n)|<C$ for some $C>0.$ 
Then we must have $(r_n+r'_n)\to 0.$
\end{prop}
\begin{proof}
First note that as we have showed $|\tr(\alpha^{k_n}_n\beta_n)|<C$ implies,
\[|z_{\alpha^{k_n}_n\beta_n,+}-z_{\alpha^{k_n}_n\beta_n,-}|\to 0\quad\text{and}\quad |z_{\alpha^{k_n}_n\beta_n,\pm}|\to 1. \]
Set $\psi_n(x)=\frac{x-z_{\alpha^{k_n}_n\beta_n,+}}{x-z_{\alpha^{k_n}_n\beta_n,-}}.$ Let $\lambda_{\alpha^{k_n}_n\beta_n}$ be the mutiplier of
$\psi_n\alpha_n^{k_n}\beta_n\psi^{-1}_n.$ By Remark \ref{4A}.1.A applied to $\psi_n<\alpha_n,\alpha^{k_n}_n\beta_n>\psi^{-1}_n$ we have,
\[|\lambda_{\alpha_n}|^2>\left(\frac{|\lambda_{\alpha^{k_n}_n\beta_n}|^{2D_n}+3}{|\lambda_{\alpha^{k_n}_n\beta_n}|^{2D_n}-1}\right)^{\frac{1}{D_n}}
\left(e^{-2\dis(\mathcal{L}_{\psi_n\alpha_n\psi_n^{-1}},\mathcal{L}_{\psi_n\alpha^{k_n}_n\beta_n\psi^{-1}_n})}\right). \] 
Since $\{z_{\psi_n\alpha_n\psi^{-1},+},z_{\psi_n\alpha_n\psi^{-1}_n,-}\}=\{1,\frac{z_{\alpha^{k_n}_n\beta_n,+}}{z_{\alpha^{k_n}_n\beta_n,-}}\}$
we have,
\[|z_{\psi_n\alpha_n\psi_n^{-1},\pm}|\to 1,\quad\text{and}\] 
\[\left|z_{\psi_n\alpha_n\psi^{-1}_n,+}-z_{\psi_n\alpha_n\psi^{-1}_n,-}\right|=\left|1-\frac{z_{\alpha^{k_n}_n\beta_n,+}}{z_{\alpha^{k_n}_n\beta_n,-}}\right|\to 0.\]
By Lemma \ref{delta-z} we have for large $n$,
\[\dis(\mathcal{L}_{\psi_n\alpha_n\psi_n^{-1}},\mathcal{L}_{\psi_n\alpha^{k_n}_n\beta_n\psi^{-1}_n})<\log\left(\frac{\delta}
{|z_{\psi_n\alpha_n\psi^{-1}_n,+}-z_{\psi_n\alpha_n\psi^{-1}_n,-}|}\right).\]
This implies that for large $n$,
\begin{eqnarray*}
|\lambda_{\alpha_n}|&>&\delta^{-1}\left(\frac{ |\lambda_{\alpha^{k_n}_n\beta_n}|^{2D_n}+3} {|\lambda_{\alpha^{k_n}_n\beta_n}|^{2D_n}-1}\right)^{\frac{1}{2D_n}}
 |z_{\psi_n\alpha_n\psi^{-1}_n,+}-z_{\psi_n\alpha_n\psi^{-1}_n,-}|\\
&>&\delta^{-1}\left(|\lambda_{\alpha^{k_n}_n\beta_n}|^{2D_n}+3\right)^{\frac{1}{2D_n}}\left(\frac{|z_{\psi_n\alpha_n\psi^{-1}_n,+}-z_{\psi_n\alpha_n\psi^{-1}_n,-}|}
{|\lambda_{\alpha^{k_n}_n\beta_n}|^2-1}\right)\\
&>&\delta^{-1}\left(|\lambda_{\alpha^{k_n}_n\beta_n}|^{2D_n}+3\right)^{\frac{1}{2D_n}}\frac{(r_n+r'_n)}
{2C'|z_{\alpha^{k_n}_n\beta_n,-}|}
\end{eqnarray*}
The last inequality follows from Proposition \ref{non-iso} and $|\lambda_{\alpha^{k_n}_n\beta_n}|<C'.$ The second inequality in the above calculations follows from that 
$|\lambda_{\alpha^{k_n}_n\beta_n}|<C'$
and for large $n$ we have, 
\[\left(|\lambda_{\alpha^{k_n}_n\beta_n}|^{2D_n}-1\right)^{\frac{1}{2D_n}}\le|\lambda_{\alpha^{k_n}_n\beta_n}|^{2D_n}-1\le|\lambda_{\alpha^{k_n}_n\beta_n}|^2-1.\]
Since $|\lambda_{\alpha_n}|<M$ for some $M,$ hence by we have,
\[r_n+r'_n<\frac{2C'M\delta|z_{\alpha^{k_n}_n\beta_n,-}|}{\left(|\lambda_{\alpha^{k_n}_n\beta_n}|^{2D_n}+3\right)^{\frac{1}{2D_n}}}<\frac{2C'M\delta'}{4^{\frac{1}{2D_n}}}\to 0.\]
\end{proof}
Now we can continue and finish the proof for $|\tr(\alpha_n^{k_n}\beta_n)|<C.$ By Proposition \ref{delta-z-r} and 
$|\lambda_{n}|>\lambda>1$ (this is condition of $(A_2)$) we have, $\frac{r_n+r'_n}{|\lambda_n|^2-1}\to 0.$ 
Since the circles $\mathcal{S}_{o_n,r_n}$ contains one of $z_{\alpha^{k_n}_n\beta_n,-},z_{\alpha^{k_n}\beta_n,+}$ and $\mathcal{S}_{o'_n,r'_n}$ contains
the other fixed point, and $\eta_{\alpha^{k_n}_n\beta_n}=1$ and $|\zeta_{\alpha^{k_n}_n\beta_n}|\to 1$ (condition of $(A_2)$), it follows from Lemma \ref{fix-trace}, we must have $\mathcal{S}_{o_n,r_n}, \mathcal{S}_{o'_n,r'_n}$
contained in the region between $\frac{1}{|\lambda_n|}$ and $|\lambda_n|$ for large $n.$ Hence we have classical generators for large $n.$
This completes the proof for $(A_2)$.\par
\end{proof}
\subsection{Case $(B)$}
Set Mobius transformations $\psi_n(x)=\zeta^{-1}_n\lambda^{2-2k_n}_n x$ and consider the generators
$\psi_n<\alpha_n,\alpha^{k_n-1}\beta_n>\psi^{-1}_n.$ Then $\zeta_{\psi_n\alpha_n^{k_n-1}\beta_n\psi_n^{-1}}=1$ and 
$\eta_{\psi_n\alpha_n^{k_n}\beta_n\psi_n^{-1}}=\zeta^{-1}_n\lambda^{2-2k_n}_n.$ Since $1\le \zeta_n\lambda_n^{2k_n}<\lambda_n^2$ and
$|\zeta_n\lambda_n^{2k_n}|-|\lambda_n^2|\to 0$ (condition of $(B)$) we have, 
\[|\zeta_n\lambda_n^{2k_n}|(1-|\zeta^{-1}_n\lambda_n^{2-2k_n}|)\to 0.\]
Since $|\lambda_n|<M$ for some $M>0,$ we have $|\eta_{\psi_n\alpha^{k_n-1}_n\beta_n\psi^{-1}_n}|\to 1.$ Hence by considering
$<\psi_n\alpha_n\psi^{-1}_n,\psi_n\beta^{-1}_n\alpha^{1-k_n}_n\psi^{-1}_n>$ we have
$|\zeta_{\psi_n\beta^{-1}_n\alpha^{1-k_n}_n\psi^{-1}_n}|\to 1$ and $\eta_{\psi_n\beta^{-1}_n\alpha^{1-k_n}_n\psi^{-1}_n}=1.$
We have reduced $(B)$ to $(A)$.\par 
Hence we have completed our proof of Theorem \ref{t-space}.
\end{proof}
\section{$\Gamma$ with small $D_\Gamma$ is either classical or there exists a universal lower bound on $Z_\Gamma$} 
This section is devoted to proving Theorem \ref{2-fixed-point} which will enable us to remove the constraint on $Z_{\G_n}$ that was placed 
in the previous section.\par
\begin{thm}\label{2-fixed-point}
There exists $c>0$ such that: Let $\Gamma_n$ be a sequence of Schottky groups with $D_n\to 0$.
Then for sufficiently large $n$, either there exists a subsequence $\Gamma_{n_i}$ that are 
classical Schottky groups, or there exists a subsequence $\Gamma_{n_j}$ with generating set 
$<\alpha_{n_j},\beta_{n_j}>$ such that $Z_{<\alpha_{n_j},\beta_{n_j}>}>c.$
\end{thm}
\begin{proof}
We prove by contradiction. Suppose there exists $\G_n$ a sequence of Schottky groups such that
for every generating set $<\alpha_n,\beta_n>$ of $\Gamma_n$ we have $Z_{<\alpha_n,\beta_n>}\to 0$.\par
For each $n$, by replacing $<\alpha_n,\beta_n>$ with $<\alpha_n,\alpha^{m_n}_n\beta_n>$ for sufficiently large $m_n$ 
if necessary, we can always assume that every $\Gamma_n$ is generated by generators with 
$|\tr(\alpha_n)|\le|\tr(\beta_n)|$ and $|\tr(\beta_n)|\ge\frac{\log3}{D_n}$. \par
We take the upper space model $\mathbb{H}^3$. By conjugating with Mobius transformations,
we can assume that $\alpha_n$ have fixed points $0,\infty$ with multiplier $\lambda_n$, and
$\beta_n$ with $z_{\beta_n,u}=1$. Recall that, as before we denote the two fixed points of
$\beta_n$ by $z_{\beta_n,l},z_{\beta_n,u}$, with $|z_{\beta_n,l}|\le|z_{\beta_n,u}|$. When we write
$\beta_n$ in matrix form, we assume that $|a_n|\le|d_n|$, otherwise we replace $\beta_n$ with $\beta^{-1}_n$.\par
By assumption, we have two cases: either case $(A)$ $z_{\beta_n,l}\to 1$ or case $(B)$ $z_{\beta_n,l}\to 0$. 
First we consider the case $(A).$\par
\subsection{Case $(A)$}
There are two possibilities, $(A_1)$ $\liminf_n|\lambda_n|=1$ or $(A_2)$ there exists
$\lambda>1$ such that $|\lambda_n|>\lambda$.
\subsubsection{$(A_2)$}
Since $|\tr(\beta_n)|\to\infty$ and $z_{\beta_n,l}\to z_{\beta_n,u}=1$ and $|\lambda_n|>\lambda>1,$ we have
$\frac{1}{\lambda}<|z_{\beta_n,l}|\le|z_{\beta_n,u}|<\lambda$ for large $n.$ Hence $<\alpha_n,\beta_n>$ satisfies Lemma \ref{classical}
for large $n.$
\subsubsection{$(A_1)$}
Taking a subsequence if necessary, we may assume that $|\lambda_n|$ is strictly decreasing
to $1$. For large enough $n$, we choose a sequence of positive integers $m_n$ depends on $n$ such that 
$1+\frac{1}{m_n+1}\le |\lambda_n|\le 1+\frac{1}{m_n}$. Let us set $\zeta_n=\frac{a_n}{c_n},\eta_n=\frac{-d_n}{c_n}$.
Since $|\tr(\beta_n)|\to\infty$ and $|\frac{\sqrt{\tr^2(\beta_n)-4}}{2c_n}|=|z_{\beta_n,+}-z_{\beta_n,-}|\to 0$
implies $|\tr(\beta_n)|<c_n,$ it follows from Lemma \ref{fix-trace} and Remark \ref{rem-2},
\[\left|z_{\alpha^{m_n}_n\beta_n,\pm}-\zeta_n\lambda^{2m_n}_n\right|\le \rho\frac{|\lambda^{m_n}_n|}{|\tr(\beta_n)|}
\le \rho\frac{e}{|\tr(\beta_n)|}\]
for large $n.$
Also by Lemma \ref{fix-trace} and Remark \ref{rem-2}, and the assumption that $z_{\beta_n,l}\to z_{\beta_n,u}=1$, we have
both $\zeta_n,\eta_n\to 1$. Hence 
\[|z_{\alpha^{m_n}_n\beta_n,\pm}-\lambda^{2m_n}_n|+|z_{\alpha^{m_n}_n\beta_n,\mp}-1|\to 0.\]
Since by our choice of $m_n$, we have $|\lambda^{2m_n}_n|\to e^2$. It follows that 
$|z_{\alpha^{m_n}_n\beta_n,\pm}|\to e^2$, and $|z_{\alpha^{m_n}_n\beta_n,\mp}|\to 1.$ Therefore, there exists $c>0$ such that
$Z_{<\alpha_n,\alpha^{m_n}_n\beta_n>}>c$ for sufficiently large $n$.
\subsection{Case $(B)$}
Here we have either: $(B_1)$ $\liminf_n|\lambda_n|< \Lambda$ for some $\Lambda>1$, or
$(B_2)$ $\liminf_n|\lambda_n|\to\infty$. We also assume that $|\zeta_n|\le|\eta_n|$ as before. 
\subsubsection{$(B_1)$} 
We will show that there exists integers $k_n$ such that $Z_{<\alpha_n,\alpha^{k_n}_n\beta_n>}>c$ for some $c>0.$\par
Take subsequence if necessary, we may assume that $|\lambda_n|\le\Lambda$ for large $n$. Choose
positive integers $k_n$ to be the smallest such that $e^2\le|\zeta_{n}\lambda^{2k_n}_n|.$ Since $|\lambda_n|\le\Lambda,$
we must have some $\sigma>0$ such that $e^2<|\zeta_{n}\lambda_n^{2k_n}|<\sigma.$
We claim that there exists $0<\epsilon<e^2$ and $n>N_\epsilon$ such that 
$e^2-\epsilon<|z_{\alpha^{k_n}_n\beta_n,+}|,|z_{\alpha^{k_n}_n\beta_n,-}|<\sigma+\epsilon$.
To see this, we use Remark \ref{4B}.2.B of Lemma \ref{fix-trace}. 
\begin{proof}{(of the claim)} 
Note that since $|z_{\beta_n,-}-z_{\beta_n,+}|<1+\epsilon$ for some $\epsilon>0$ and large $n,$ and also $|\tr(\beta_n)|\to\infty,$ we have
$|c_n|\to\infty.$ By Remark \ref{4B}.2.B we have $|z_{\beta_n,u}-\eta_n|\to 0,$ hence $\eta_n\to 1.$\par
First we show that $|z_{\alpha^{k_n}_n\beta_n,+}-z_{\alpha^{k_n}\beta_n,-}|\not\to\infty.$ Assume otherwise.
Let $\rho_n$ be the center of the circle having 
$z_{\alpha^{k_n}_n\beta_n,+}$ and $z_{\alpha^{k_n}_n\beta_n,-}$ as antipodal points. Since 
$|\zeta_{\alpha^{k_n}_n\beta_n}+\eta_{\alpha^{k_n}_n\beta_n}|<\sigma+1+\epsilon'$ for some $\epsilon'>0$ and large $n,$ 
and $\rho_n=\frac{z_{\alpha^{k_n}_n\beta_n,+}+z_{\alpha^{k_n}\beta_n,-}}{2},$ and 
$z_{\alpha^{k_n}_n\beta_n,+}+z_{\alpha^{k_n}\beta_n,-}=\zeta_{\alpha^{k_n}_n\beta_n}+\eta_{\alpha^{k_n}_n\beta_n}$ we have 
$|\rho_n-\sigma'|<\kappa$ for some $\kappa>0$ and $\sigma'=\sigma+1+\epsilon'.$ Note that since
$|\rho_n-\sigma'|\le|\rho_n|+|\sigma'|<\frac{\sigma+1+\epsilon'}{2}+\sigma+1+\epsilon'$ we can take $\kappa=\frac{3}{2}(\sigma+1+\epsilon').$
This
imples that $\dis(\mathcal{L}_{\alpha_n},\mathcal{L}_{\alpha^{k_n}_n\beta_n})<\delta$ for some $\delta>0.$ By Remark \ref{4A}.1.A we have 
$|\tr(\alpha^{k_n}_n\beta_n)|\to\infty,$ 
and $|\zeta_{\alpha^{k_n}_n\beta_n}-\eta_{\alpha^{k_n}_n\beta_n}|=\frac{|\tr(\alpha^{k_n}_n\beta_n)|}{|c_n\lambda^{-k_n}_n|},$ we get
$|\tr(\alpha^{k_n}_n\beta_n)|\asymp|c_n\lambda^{-k_n}_n|.$ But this imples that 
$|z_{\alpha^{k_n}_n\beta_n,+}-z_{\alpha^{k_n}\beta_n,-}|<C$ for some $C>0,$ hence a contradiction.\par
Note that if $|z_{\alpha^{k_n}\beta_n,+}-z_{\alpha^{k_n}\beta_n,-}|\to 0$ then
$|z_{\alpha^{k_n}\beta_n,\pm}|\to \frac{1}{2}|\zeta_{\alpha^{k_n}_n\beta_n}+\eta_{\alpha_n^{k_n}\beta_n}|
=\frac{1}{2}|\zeta_{\alpha^{k_n}_n\beta_n}+1|.$ Since $\frac{1}{2}(e^2-1)<\frac{1}{2}|\zeta_{\alpha^{k_n}_n\beta_n}+1|<\frac{1}{2}\sigma'$ for 
large $n.$  
This implies that $\frac{1}{2}(e^2-1)<|z_{\alpha^{k_n}\beta_n,\pm}|<\frac{1}{2}\sigma'$ for large $n.$\par
Finally if $c<|z_{\alpha^{k_n}_n\beta_n,+}-z_{\alpha^{k_n}\beta_n,-}|<c'$ for some $c,c'>0,$ then by Remark \ref{4A}.1.A we have 
$|\tr(\alpha^{k_n}_n\beta_n)|\to\infty$ which imples $|c_n\lambda_n^{-k_n}|\to\infty.$ Hence by Remark \ref{4B}.2.B we have
$\{z_{\alpha^{k_n}_n\beta_n,+},z_{\alpha^{k_n}\beta_n,-}\}\to\{\zeta_{\alpha^{k_n}_n\beta_n},\eta_{\alpha_n^{k_n}\beta_n}\}$ which implies the claim.
\end{proof}
Now with the claim true, there are two possibilities: $(B'_1)$  
$\liminf|z_{\alpha^{k_n}_n\beta_n,+}-z_{\alpha^{k_n}_n\beta_n,-}|\to 0$, or $(B''_1)$
$\liminf_n|z_{\alpha^{k_n}_n\beta_n,+}-z_{\alpha^{k_n}_n\beta_n,-}|>0$. 
For $(B''_1)$, we have $Z_{<\alpha_n,\alpha^{k_n}_n\beta_n>}>c$, for some $c>0$ and large $n$.\par
Suppose $(B'_1).$ By passing to subsequence if necessary, we take $|z_{\alpha^{k_n}_n\beta_n,+}-z_{\alpha^{k_n}_n\beta_n,-}|\to 0.$\par
If $|\lambda_n|\to 1$ then we choose positive integers $m_n$ as defined in case $(A_1)$. 
Then $e^2\le\zeta_n\lambda^{2k_n+2m_n}_n<e^4,$ and $\kappa_1<|z_{\alpha^{k_n+m_n}_n\beta_n,+}-z_{\alpha^{k_n+m_n}_n\beta_n,-}|<\kappa_2$  
for some $0<\kappa_1,\kappa_2.$ By Remark \ref{4A}.1.A we have $|\tr(\alpha^{k_n+m_n}_n\beta_n)|\to\infty.$ Hence by
Remark \ref{4B}.2.B we have 
$\{z_{\alpha^{k_n+m_n}_n\beta_n,+},z_{\alpha^{k_n+m_n}_n\beta_n,-}\}\to\{\zeta_{\alpha^{k_n+m_n}_n\beta_n},\eta_{\alpha^{k_n+m_n}_n\beta_n}\}.$
This implies that $e^2\le|z_{\alpha^{k_n+m_n}_n\beta_n,\pm}|<e^4$ and $|z_{\alpha^{k_n+m_n}_n\beta_n,\mp}|\to 1$. 
Hence there exists $c>0$, such that $Z_{<\alpha_n,\alpha^{k_n+m_n}_n\beta_n>}>c$ for large $n$.\par
If $\lambda_n|>c>0$ then we take $m_n=1.$ Then 
$\delta_1<|z_{\alpha^{k_n+1}_n\beta_n,+}-z_{\alpha^{k_n+1}_n\beta_n,-}|<\delta_2$  
for some $0<\delta_1,\delta_2.$ And it follows from Remark \ref{4A}.1.A and Remark \ref{4B}.2.B we have that
$Z_{<\alpha_n,\alpha^{k_n+1}_n\beta_n>}>c$ for large $n$.
\subsubsection{$(B_2)$} 
By taking a subsequence of $\alpha_n$, we may assume that $|\lambda_n|$ is strictly increasing.
Choose a sequence of largest integers $k_n\ge 0$, such that $|\zeta_n\lambda^{2k_n}_n|\le 1$.
If $\limsup |\zeta_n\lambda^{2k_n}_n|=1$ but $\limsup|\zeta_n\lambda^{2k_n}_n-1|\not=0$, then there is a subsequence
$<\alpha_{n_j},\alpha^{k_{n_j}}_{n_j}\beta_{n_j}>$ of $<\alpha_n,\alpha^{k_n}_n\beta_n>$, such that 
$\liminf_j Z_{<\alpha_{n_j},\alpha^{k_{n_j}}_{n_j}\beta_{n_j}>}>0$.\par
If $\limsup \zeta_n\lambda^{2k_n}_n=1$, then let $<\alpha_{n_i},\alpha^{k_{n_i}}_{n_i}\beta_{n_i}>$ be the subsequence
of $<\alpha_n,\alpha^{k_n}_n\beta_n>$ with $\lim_i\zeta_{n_i}\lambda^{2k_{n_i}}_{n_i}=1$. \par
If $\limsup|\tr(\alpha^{k_{n_i}}_{n_i}\beta_{n_i})|=\infty$, then by passing to a subsequence if necessary, for large
$i$, $\alpha^{k_{n_i}}_{n_i}\beta_{n_i}$ will have disjoint isometric circles. 
Note $|\tr(\beta_{n_i})|\to\infty$, $z_{\beta_{n_i},u}=1$ and $z_{\beta_{n_i},l}\to 0,$ so $|z_{\beta_{n_i},+}-z_{\beta_{n_i},-}|\le 1+\epsilon_n$
for $\epsilon_n\to 0$ (i.e. $|z_{\beta_{n_i},+}-z_{\beta_{n_i},-}|<c,$ for $c>0$), we have $|c_{n_i}|\to\infty.$
By Remark \ref{4B}.2.B, $\lim_i\min\{|\eta_{n_i}-z_{\beta_{n_i},-}|,|\eta_{n_i}-z_{\beta_{n_i},+}|\}\to 0$. 
Since $\zeta_{\alpha^{k_{n_i}}_{n_i}\beta_{n_i}}\to 1=\eta_{\alpha^{k_{n_i}}_{n_i}\beta_{n_i}}$ and $|\tr(\alpha^{k_{n_i}}_{n_i}\beta_{n_i})|\to\infty$ we have,
$|\zeta_{\alpha^{k_{n_i}}_{n_i}\beta_{n_i}}-\eta_{\alpha^{k_{n_i}}_{n_i}\beta_{n_i}}|=\frac{|\tr(\alpha^{k_{n_i}}_{n_i}\beta_{n_i})|}
{|\lambda^{-k_{n_i}}_{n_i}c_{n_i}|}\to 0\quad\text{and}$
$|\lambda^{-k_{n_i}}_{n_i}c_{n_i}|\to\infty.$ Also 
$|z_{\alpha^{k_{n_i}}_{n_i}\beta_{n_i},+}-z_{\alpha^{k_{n_i}}_{n_i}\beta_{n_i},-}|=\frac{|\sqrt{\tr^2(\alpha^{k_{n_i}}_{n_i}\beta_{n_i})-4}|}
{2|\lambda^{-k_{n_i}}_{n_i}c_{n_i}|}\to 0.$
Hence by Remark \ref{4B}.2.B, there exists a $\kappa>0$ such that for large $i$ we have,
 $\kappa^{-1}<|z_{\alpha^{k_{n_i}}_{n_i}\beta_{n_i},l}|\le|z_{\alpha^{k_{n_i}}_{n_i}\beta_{n_i},u}|<\kappa$.
And since $|\lambda_{n_i}|\to\infty$,
we can choose Mobius transformations $\psi_i$ such that 
$\psi_i<\alpha_{n_i},\alpha^{k_{n_i}}_{n_i}\beta_{n_i}>\psi^{-1}_i$ satisfies Lemma \ref{classical}. \par
If $\limsup|\tr(\alpha^{k_{n_i}}_{n_i}\beta_{n_i})|<\infty$, then let $\phi_i$ be Mobius transformations
such that $\phi_i\alpha^{k_{n_i}}_{n_i}\beta_{n_i}\phi^{-1}_{i}$ have fixed points $0,\infty$, and
fixed point $z_{\phi_{i}\alpha_{n_i}\phi^{-1}_i,u}$ of $\alpha_{n_i}$ is $1$. Then it follows that 
$z_{\phi_i\alpha_{n_i}\phi^{-1}_i,l}\to 1$. Since $|\tr(\phi_i\alpha_{n_i}\phi^{-1}_i)|\to\infty$, 
hence we have reduced this case to case $(A)$, which we already considered.\par
If $\limsup|\zeta_n\lambda^{2k_n}_n|<1$, then we have two possibilities:
$(B'_2)$ $\liminf|\zeta_n\lambda^{2k_n+2}_n|=1$, or $(B''_2)$ $\liminf|\zeta_n\lambda^{2k_n+2}_n|>1.$\par
Consider case $(B'_2)$. Let $<\alpha_{n_i},\alpha^{k_{n_i}}_{n_i}\beta_{n_i}>$ be a subsequence of
$<\alpha_n,\alpha^{k_n}_n\beta_n>$ such that $\lim_i|\zeta_{n_i}\lambda^{2k_{n_i}+2}_{n_i}|\to 1$. 
If $\sup_i|\tr(\alpha^{k_{n_i+1}}_{n_i}\beta_{n_i})|<\infty$, then we conjugate 
$\alpha_{n_i},\alpha^{k_{n_i}}_{n_i}\beta_{n_i}$ to
$\hat{\beta}_i=\phi_i\alpha_{n_i}\phi^{-1}_i,\hat{\alpha}_i=\phi_i\alpha^{k_{n_i}+1}_{n_i}\beta_{n_i}\phi^{-1}_i$ with $\hat{\alpha}_i$ have fixed points $0,\infty$ and $\hat{\beta}_i$ have $z_{\hat{\beta}_i,u}=1$. 
Since $\sup_i|\hat{\lambda}_i|<\infty$, it follows that, if $z_{\hat{\beta}_i,l}\to 1$ then 
$<\hat{\alpha}_i,\hat{\beta}_i>$, falls under case $(A),$  and if
$z_{\hat{\beta}_i,l}\to 0$ then $<\hat{\alpha}_i,\hat{\beta}_i>$ falls under case $(B_1).$ Otherwise there exists $\epsilon>0$ such that
$\epsilon<|z_{\hat{\beta}_i,l}|<1-\epsilon,$ hence $Z_{<\hat{\alpha}_i,\hat{\beta}_i>}>c$ for some $c>0.$\par
On the other hand, if $\sup_i|\tr(\alpha^{k_{n_i}+1}_{n_i}\beta_{n_i})|=\infty$, then for large $i$
since the radius of isometric circles is
$\mathfrak{R}_{\alpha^{k_{n_i}+1}_{n_i}\beta_{n_i}}=\frac{\left|z_{\alpha^{k_{n_i}+1}_{n_i}\beta_{n_i},u}-z_{\alpha^{k_{n_i}+1}_{n_i}\beta_{n_i},l}\right|}
{|\tr(\alpha^{k_{n_i}+1}_{n_i}\beta_{n_i})|}$ and the distance between the centers of these isometric circles is 
$|\zeta_{\alpha^{k_{n_i}+1}_{n_i}\beta_{n_i}}-\eta_{\alpha^{k_{n_i}+1}_{n_i}\beta_{n_i}}|=$
$\frac{|\tr(\alpha^{k_{n_i}+1}_{n_i}\beta_{n_i})|}{c_n\lambda_{n_i}^{-k_{n_i}-1}},$ and by
$\left|z_{\alpha^{k_{n_i}+1}_{n_i}\beta_{n_i},u}-z_{\alpha^{k_{n_i}+1}_{n_i}\beta_{n_i},l}\right|=
\frac{\sqrt{\tr^2(\alpha^{k_{n_i}+1}_{n_i}\beta_{n_i})-4}}{2c_{n_i}\lambda_{n_i}^{-k_{n_i}-1}}$
we have,
$\lim_i\frac{\left|\zeta_{\alpha^{k_{n_i}+1}_{n_i}\beta_{n_i}}-\eta_{\alpha^{k_{n_i}+1}_{n_i}\beta_{n_i}}\right|}{\mathfrak{R}_{\alpha^{k_{n_i}+1}_{n_i}\beta_{n_i}}}
=\lim_i\frac{2\left|\tr(\alpha^{k_{n_i}+1}_{n_i}\beta_{n_i})\right|^2}{\sqrt{\tr^2(\alpha^{k_{n_i}+1}_{n_i}\beta_{n_i})-4}}>\delta
|\tr(\alpha^{k_{n_i}+1}_{n_i}\beta_{n_i})|$ for some $\delta>0.$ Hence
$\left|\zeta_{\alpha^{k_{n_i}+1}_{n_i}\beta_{n_i}}-\eta_{\alpha^{k_{n_i}+1}_{n_i}\beta_{n_i}}\right|>2\mathfrak{R}_{\alpha^{k_{n_i}+1}_{n_i}\beta_{n_i}}$
for large $i.$
This implies $\alpha^{k_{n_i}+1}_{n_i}\beta_{n_i}$ have disjointed isometric circles for large $i.$ 
By Lemma \ref{fix-trace} and $z_{\beta_{n_i},u}=1$ we have  $\eta_{\beta_{n_i}}\to 1.$ And since $\eta_{\alpha^{k_{n_i}+1}_{n_i}\beta_{n_i}}=
\eta_{\beta_{n_i}},$ implies that $\eta_{\alpha^{k_{n_i}+1}_{n_i}\beta_{n_i}}\to 1.$ 
Note that if $\inf_i|\zeta_{\alpha^{k_{n_i}+1}_{n_i}\beta_{n_i}}-1|>0$ and $|\zeta_{\alpha^{k_{n_i}+1}_{n_i}\beta_{n_i}}|\to 1,$ then
by $|\tr(\alpha^{k_{n_i}+1}_{n_i}\beta_{n_i})|\to\infty$ we have 
$\left|1-z_{\alpha^{k_{n_i}+1}_{n_i}\beta_{n_i},l}\right|,$ 
$\left|\zeta_{\alpha^{k_{n_i}+1}_{n_i}\beta_{n_i}}-z_{\alpha^{k_{n_i}+1}_{n_i}\beta_{n_i},u}\right|\to 0,$ and 
$\inf_i\left|z_{\alpha^{k_{n_i}+1}_{n_i}\beta_{n_i},u}-z_{\alpha^{k_{n_i}+1}_{n_i}\beta_{n_i},l}\right|>0.$ Hence for large $i$ there exits
$\epsilon>0$ such that $1-\epsilon<|z_{\alpha^{k_{n_i}+1}_{n_i}\beta_{n_i},l}|\le|z_{\alpha^{k_{n_i}+1}_{n_i}\beta_{n_i},u}|<1+\epsilon,$
therefore we must have $\inf_iZ_{<\alpha_{n_i},\alpha^{k_{n_i}+1}_{n_i}\beta_{n_i}>}>0.$
It follows we can assume that $\zeta_{\alpha^{k_{n_i}+1}_{n_i}\beta_{n_i}}\to 1.$ Then
we have $|\zeta_{\alpha^{k_{n_i}+1}_{n_i}\beta_{n_i}}-\eta_{\alpha^{k_{n_i}+1}_{n_i}\beta_{n_i}}|\to 0.$ Hence 
$\frac{|\tr(\alpha^{k_{n_i}+1}_{n_i}\beta_{n_i})|}{c_n\lambda_{n_i}^{-k_{n_i}-1}}\to 0.$ Since 
$|\tr(\alpha^{k_{n_i}+1}_{n_i}\beta_{n_i})|\to\infty,$ which gives
$\left|z_{\alpha^{k_{n_i}+1}_{n_i}\beta_{n_i},u}-z_{\alpha^{k_{n_i}+1}_{n_i}\beta_{n_i},l}\right|=
\frac{\sqrt{\tr^2(\alpha^{k_{n_i}+1}_{n_i}\beta_{n_i})-4}}{2c_{n_i}\lambda_{n_i}^{-k_{n_i}-1}}$ $\to 0.$
Therefore the distance between the centers of these isometric circles decreases to $0$ and radius $\mathfrak{R}_{\alpha^{k_{n_i}+1}_{n_i}\beta_{n_i}}\to 0.$
Since $\inf_i|\lambda_{n_i}|^2>c>1,$ we have for large $i$ that the isometric circles of $\alpha^{k_{n_i}+1}_{n_i}\beta_{n_i}$ disjointed and
lies between $c^{-1}$ and $c.$ In particular, $c^{-1}<|z_{\alpha^{k_{n_i}+1}_{n_i}\beta_{n_i},l}|\le|z_{\alpha^{k_{n_i}+1}_{n_i}\beta_{n_i},u}|
<c,$ and satisfies Lemma \ref{classical}.\par
Case $(B''_2)$. First we define a new sequence of $<\tilde{\alpha}_{n},\tilde{\beta}_{n}>$ as follows:
Consider $<\alpha_n,\alpha^{k_n}_n\beta_n>$. If $|\tr(\alpha^{k_n}_n\beta_n)|\ge|\tr(\alpha_n)|$,
then set $\tilde{\alpha}_n=\alpha_n, \tilde{\beta}_n=\beta_n$. Otherwise, let $\phi_n$ be the Mobius map 
so that $\phi_n\alpha^{k_n}_n\beta_n\phi^{-1}_n$ 
have fixed points $0,\infty$, and $\phi_n\alpha_n\phi^{-1}_n$ have $z_{\phi_n\alpha_n\phi^{-1}_n,u}=1$.
Set $\alpha_{n,1}=\phi_n\alpha^{k_n}_n\beta_n\phi^{-1}_n, \beta_{n,1}=\phi_n\alpha_n\phi^{-1}_n$.
We define integer $k_{n,1}$ with respect to $<\alpha_{n,1},\beta_{n,1}>$ the same way as we defined $k_n$ before. \par
Now if $|\tr(\alpha^{k_{n,1}}_{n,1}\beta_{n,1})|\ge|\tr(\alpha_{n,1})|$ then 
we set $\tilde{\beta}_n=\beta_{n,1}$ and $\tilde{\alpha}_n=\alpha_{n,1}$. Otherwise, we repeat this construction
to get a sequence $<\alpha_{n,m},\beta_{n,m}>$. 
By construction for a each $n$, either there exists a $m$ such that $|\tr(\alpha^{k_{n,m}}_{n,m}\beta_{n,m})|\ge|\tr(\alpha_{n,m})|$ 
or we have $|\tr(\alpha^{k_{n,m}}_{n,m}\beta_{n,m})|<|\tr(\alpha_{n,m})|$ for all $m.$ Assume the latter holds, since
$\alpha_{n,m+1}=\phi_{n,m}\alpha^{k_{n,m}}_{n,m}\beta_{n,m}\phi^{-1}_{n,m}$ we have 
$|\tr(\alpha_{n,m+1})|<|\tr(\alpha_{n,m})|$ for all $m.$ If $\lim_m|\tr(\alpha_{n,m})|=0$ then take $m_n$ to be the first integer $m$ with
with $|\tr(\alpha_{n,m})|<\frac{1}{n}.$ If $\lim_m|\tr(\alpha_{n,m})|>0$ then take $m_n$ to be the first integer $m$ with 
$|\tr(\alpha_{n,m+1})|>|\tr(\alpha_{n,m})|-\frac{1}{n}.$ If the former holds, we set $m_n$ to be the first integer $m$ with
$|\tr(\alpha^{k_{n,m}}_{n,m}\beta_{n,m})|\ge|\tr(\alpha_{n,m})|.$ Hence there exists a $m_n$, such that either
$|\tr(\alpha^{k_{n,m_n}}_{n,m_n}\beta_{n,m_n})|>|\tr(\alpha_{n,m_n})|-\frac{1}{n}$, or 
$|\tr(\alpha_{n,m_n})|<\frac{1}{n}$. We define $\tilde{\alpha}_n=\alpha_{n,m_n},\tilde{\beta}_{n}=\beta_{n,m_n}.$\par

Now consider $<\tilde{\alpha}_n,\tilde{\beta}_n>$. If $\liminf_n|\tr(\tilde{\alpha}_n)|<\infty$, we choose a subsequence
with $|\tr(\tilde{\alpha}_{n_i})|<c$ for all large $i$ and some $c>0$. Let $p_i$ be a sequence of least positive integers 
such that $|\tr(\tilde{\alpha}^{p_i}_{n_i}\tilde{\beta}_{n_i})|>\frac{1}{D_{n_i}}$. We conjugate $\tilde{\alpha}^{p_i}_{n_i}\tilde{\beta}_{n_i}$
by $\psi_i$ that fixes $0,\infty$ and $z_{\psi_i\tilde{\alpha}^{p_i}_{n_i}\tilde{\beta}_{n_i}\psi^{-1}_i,u}=1$. 
Set $\bar{\alpha}_i=\psi_i\tilde{\alpha}_{n_i}\psi^{-1}_i, \bar{\beta}_i=\psi_i\tilde{\alpha}^{p_i}_{n_i}\tilde{\beta}_{n_i}\psi^{-1}_i.$
By construction, if $z_{\bar{\alpha}_i,l}\to 0$ then $<\bar{\alpha}_i,\bar{\beta}_i>$  satisfies $(B_1),$ and if
$z_{\bar{\alpha}_i,l}\to 1$ then $<\bar{\alpha}_i,\bar{\beta}_i>$  satisfies $(A).$ Otherwise there exists $\epsilon>0$
such that $\epsilon<|z_{\bar{\alpha}_i,l}|<1-\epsilon$ which implies that $Z_{<\bar{\alpha}_i,\bar{\beta}_i>}>c$ for some $c>0.$
Hence in either case, we are done. \par
On the other hand, suppose $\liminf_n|\tr(\tilde{\alpha}_n)|=\infty$. Since $|\tr(\tilde{\beta}_n)|\ge|\tr(\tilde{\alpha}_n)|$ then
it's sufficient to assume that $<\tilde{\alpha}_n,\tilde{\beta}_n>$ satisfies case $(B''_2)$, otherwise
we are done. We define $\tilde{k}_n=k_{n,m_n}$. \par
Set $\nu_n=\tilde{\alpha}^{\tilde{k}_{n}}_{n}\tilde{\beta}_{n},\mu_n=\tilde{\alpha}_{n}.$
Since $|\tr(\tilde{\beta}_n)|\ge|\tr(\tilde{\alpha}_n)|$ and 
$|z_{\tilde{\beta}_n,-}-z_{\tilde{\beta}_n,+}|\le 1+\delta_n$ with $\delta_n\to 0$ (this follows from that 
$z_{\tilde{\beta}_n,u}=1, z_{\tilde{\beta}_n,l}\to 0$), 
it follows from Lemma \ref{fix-trace} with Remark \ref{rem-2} and $|\tr(\alpha^{k_{n,m_n}}_{n,m_n}\beta_{n,m_n})|>|\tr(\alpha_{n,m_n})|-\frac{1}{n}$
which implies $\lim_n\frac{|\tr(\nu_n)|}{|\tilde{\lambda}_{n}|}>\epsilon$ for $\epsilon>0$, we have
$|\eta_{\nu_n}-1|<\delta|\tilde{\lambda}_{n}|^{-2},$ for some $\delta>0.$ 
Since $\eta_{\nu_n\mu_n}=\eta_{\nu_n}\tilde{\lambda}_n^{-2}$ and $\eta_{\mu_n\nu_n}=\eta_{\nu_n},$ we have
$|\eta_{\nu_n\mu_n}-\tilde{\lambda}_n^{-2}|=|\eta_{\nu_n}\tilde{\lambda}_n^{-2}-\tilde{\lambda}_n^{-2}|<\delta|\tilde{\lambda}_{n_n}|^{-4}$ and
$|\eta_{\mu_n\nu_n}-1|=|\eta_{\nu_n}-1|<\delta|\tilde{\lambda}_{n_n}|^{-2},$ for large $n$\par
We have two cases to consider:
\[(1):\quad\quad\lim_n\frac{|\tilde{\zeta}_n\tilde{\lambda}^{2\tilde{k}_n}_n|-|\tilde{\lambda}_n|^{-2}}{|\tilde{\lambda}_n|^{-2}}=\infty.
\quad\hspace{4.3cm}\]
\[(2):\quad\quad\lim_i\frac{|\tilde{\zeta}_{n_i}\tilde{\lambda}^{2\tilde{k}_{n_i}}_{n_i}|-|\tilde{\lambda}^{-2}_{n_i}|}{|\tilde
{\lambda}_{n_i}|^{-2}}<\sigma,\quad\text{for some subsequence $n_i$.}\]
\begin{proof}{(Assume (2))}
Since $|\tilde{\zeta}_{n_i}\tilde{\lambda}^{2\tilde{k}_{n_i}}_{n_i}|>|\tilde{\lambda}_{n_i}|^{-2}$ and $(2),$
and $\zeta_{\mu_i\nu_i}=\tilde{\lambda}_{n_i}^2\zeta_{\nu_i},$ we have for large $i$,
\[1<|\zeta_{\mu_i\nu_i}|<(\sigma+1).\]
Since we also have $1-\delta|\tilde{\lambda}_{n_i}^{-2}|<|\eta_{\mu_i\nu_i}|<1+\delta|\tilde{\lambda}_{n_i}^{-2}|$ and 
$|\tilde{\lambda}_{n_i}|\to\infty$ it follows that
there exists $\kappa>1$ such that, $1-\kappa^{-1}<|\eta_{\mu_i\nu_i}|,|\zeta_{\mu_i\nu_i}|<1+\kappa$ for large $i.$ By Lemma \ref{fix-trace} 
with Remark \ref{rem-2} and
$|\tr(\nu_i)|>\epsilon|\tilde{\lambda}_{n_i}|$ for large $i$ we have for some $\rho_1,\rho_2>0$ that,
\[1-\kappa^{-1}-\rho_1|\tilde{\lambda}_{n_i}|^{-1}<|z_{\mu_i\nu_i,l}|,|z_{\mu_i\nu_i,u}|<1+\kappa+\rho_2|\tilde{\lambda}_{n_i}|^{-1}.\]
Hence there exists $\kappa'>1$ such that $\kappa'^{-1}<|z_{\mu_i\nu_i,l}|\le|z_{\mu_i\nu_i,u}|<\kappa'.$\par
If $|\tr(\mu_i\nu_i)|\to\infty$, then $<\mu_i,\mu_i\nu_i>$ satisfies the second set of conditions of Lemma \ref{classical},
hence classical.\par
If $\limsup|\tr(\mu_i\nu_i)|<\infty$, then define Mobius transformations $\psi_i$ such that
$\psi_i\mu_i\nu_i\psi^{-1}$ have fixed points $0,\infty$ and $z_{\psi_i\mu_i\psi_i^{-1},u}=1.$
If $z_{\psi_i\mu_i\psi_i^{-1},l}\to 0$ or $z_{\psi_i\mu_i\psi_i^{-1},l}\to 1$ then 
$\psi_i<\mu_i,\mu_i\nu_i>\psi^{-1}_i$ satisfies $(A)$ or $(B_1).$ Otherwise we have for some $\epsilon>0$ such that,
$\epsilon<|z_{\psi_i\mu_i\psi^{-1},l}|<1-\epsilon,$ which implies that $Z_{\psi_i<\mu_i,\mu_i\nu_i>\psi^{-1}_i}>c$
for some $c>0.$ This completes our proof of $(B_2'')$ with $(2)$.
\end{proof}
\begin{proof}{(Assume $(1)$)}
By $(1),$ there exists $0<\rho_n\to\infty$ with $\rho_n<|\tilde{\lambda}_n|$, such that 
$|\tilde{\zeta}_n\tilde{\lambda}^{2\tilde{k}_n}_n|-|\tilde{\lambda}_n|^{-2}>\rho_n|\tilde{\lambda}_n|^{-2}.$ 
Let $\chi_n$ be Mobius transformations defined by $\chi_n(x)=\frac{\tilde{\lambda}_n}{\sqrt{\rho_n}}x.$ 
We will show that $\chi_n<\mu_n,\nu_n>\chi_n^{-1}$ satisfies Remark \ref{classical-rem-2} of Lemma \ref{classical}.  \par
Since $|\eta_{\nu_n}-1|<\frac{\delta}{|\tilde{\lambda}_n|^2}$ we have,
\[\frac{|\tilde{\lambda}_n|}{\sqrt{\rho_n}}-\frac{\delta}{|\tilde{\lambda}_n|\sqrt{\rho_n}}<|\eta_{\chi_n\nu_n\chi_n^{-1}}|<
\frac{|\tilde{\lambda}_n|}{\sqrt{\rho_n}}+\frac{\delta}{|\tilde{\lambda}_n|\sqrt{\rho_n}}.\]
By condition of $(B_2)$ we have $|\tilde{\zeta}_n\tilde{\lambda}_n^{2\tilde{k}_n}|<1.$ This gives,
\[\frac{\sqrt{\rho_n}}{|\tilde{\lambda}_n|}+\frac{1}{|\tilde{\lambda}_n|\sqrt{\rho_n}}<|\zeta_{\chi_n\nu_n\chi_n^{-1}}|<\frac{|\tilde{\lambda}_n|}{\sqrt{\rho_n}}.\]

By Lemma \ref{fix-trace} with Remark \ref{rem-2}, $|z_{\nu_i,\pm}-\zeta_{\nu_n}|<\sigma|\tr(\nu_n)|^{-2}$ and
$|z_{\nu_i,\mp}-\eta_{\nu_n}|<\sigma|\tr(\nu_n)|^{-2}$ for some $\sigma>0.$ Since $|\tr(\nu_n)|>\epsilon|\tilde{\lambda}_n|$ for large $n$
we have, 
\[|z_{\chi_n\nu_n\chi_n^{-1},\pm}-\eta_{\chi_n\nu_n\chi_n^{-1}}|<\frac{\sigma}{\epsilon|\tilde{\lambda}_n|\sqrt{\rho_n}},\quad
|z_{\chi_n\nu_n\chi_n^{-1},\mp}-\zeta_{\chi_n\nu_n\chi_n^{-1}}|<\frac{\sigma}{\epsilon|\tilde{\lambda}_n|\sqrt{\rho_n}}.\]
Hence, 
\[|z_{\chi_n\nu_n\chi_n^{-1},u}|<\frac{|\tilde{\lambda}_n|}{\sqrt{\rho_n}}+\frac{\delta}{|\tilde{\lambda}_n|\sqrt{\rho_n}}
+\frac{\sigma}{\epsilon|\tilde{\lambda}_n|\sqrt{\rho_n}}\quad\text{and},\]
\[|z_{\chi_n\nu_n\chi_n^{-1},_l}|>\frac{\sqrt{\rho_n}}{|\tilde{\lambda}_n|}-\frac{\delta}{|\tilde{\lambda}_n|\sqrt{\rho_n}}-
\frac{\sigma}{\epsilon|\tilde{\lambda}_n|\sqrt{\rho_n}}.\]
We have $|\tilde{\lambda}_n|^{-1}<|z_{\chi_n\nu_n\chi_n^{-1},l}|\le|z_{\chi_n\nu_n\chi^{-1}_n,u}|<|\tilde{\lambda}_n|$ for large $n.$\par
By above estimates for fixed points of $\chi_n\nu_n\chi_n^{-1},$ and $|\tr(\nu_n)|>\epsilon|\tilde{\lambda}_n|$ we have,
\begin{eqnarray*}
\frac{(|z_{\chi_n\nu_n\chi_n^{-1},u}|+1)(|\tilde{\lambda}_n|+1)}{|\tr(\nu_n)|(|\tilde{\lambda}_n|-|z_{\chi_n\nu_n\chi_n^{-1},u}|)}
&<&\frac{\frac{|\tilde{\lambda}_n|^2}{\sqrt{\rho_n}} +\frac{|\tilde{\lambda}_n|}{\sqrt{\rho_n}}+|\tilde{\lambda}_n|+\delta'
}{\epsilon|\tilde{\lambda}_n|^2(1-\frac{1}{\rho_n}-\frac{\delta}{|\tilde{\lambda}_n|^2\sqrt{\rho_n}}
-\frac{\sigma}{\epsilon|\tilde{\lambda}_n|^2\sqrt{\rho_n}})}\\
&<&\frac{1 +\frac{1}{|\tilde{\lambda}_n|}+\frac{\sqrt{\rho_n}}{|\tilde{\lambda}_n|}+\frac{\delta'\sqrt{\rho_n}}{|\tilde{\lambda}_n|^2}
}{\epsilon\sqrt{\rho_n}(1-\frac{1}{\rho_n}-\frac{\delta}{|\tilde{\lambda}_n|^2\sqrt{\rho_n}}
-\frac{\sigma}{\epsilon|\tilde{\lambda}_n|^2\sqrt{\rho_n}})}\\
\text{by $\rho_n<|\tilde{\lambda}_n|$} &\text{we}& \text{have some $\delta''>1$ such that,}\\
&<&\frac{\delta''}{\epsilon\sqrt{\rho_n}(1-\frac{1}{\rho_n}-\frac{\delta}{|\tilde{\lambda}_n|^2\sqrt{\rho_n}}
-\frac{\sigma}{\epsilon|\tilde{\lambda}_n|^2\sqrt{\rho_n}})}\to\ 0.
\end{eqnarray*}
For the other part of the conditions of Remark \ref{classical-rem-2} we have: \\
If $|z_{\chi_n\nu_n\chi_n^{-1},l}|<M$ then, 
\begin{eqnarray*}
\frac{(|z_{\chi_n\nu_n\chi_n^{-1},l}|+1)(|\tilde{\lambda}_n|^{-1}+1)}{|\tr(\nu_n)|(|z_{\chi_n\nu_n\chi_n^{-1},l}|-|\tilde{\lambda}_n|^{-1})}&<&
\frac{(M+1)(|\tilde{\lambda}_n|^{-1}+1)}
{\epsilon|\tilde{\lambda}_n|(\frac{\sqrt{\rho_n}}{|\tilde{\lambda}_n|}-\frac{\delta}{|\tilde{\lambda}_n|\sqrt{\rho_n}}-
\frac{\sigma}{\epsilon|\tilde{\lambda}_n|\sqrt{\rho_n}}-\frac{1}{|\tilde{\lambda}_n|})}\\
&<&\frac{M'}{\epsilon(\sqrt{\rho_n}-\frac{\delta}{\sqrt{\rho_n}}-\frac{\sigma}{\epsilon\sqrt{\rho_n}}-1)}
<\frac{M'}{\epsilon''\sqrt{\rho_n}}\to 0.
\end{eqnarray*}
Otherwise we have $|z_{\chi_n\nu_n\chi_n^{-1},l}|\to\infty$ and,
\[
\frac{(|z_{\chi_n\nu_n\chi_n^{-1},l}|+1)(|\tilde{\lambda}_n|^{-1}+1)}{|\tr(\nu_n)|(|z_{\chi_n\nu_n\chi_n^{-1},l}|-|\tilde{\lambda}_n|^{-1})}
<\frac{\delta''}{\epsilon|\tilde{\lambda}_n|}\to 0,\quad\text{for some $\delta''>0$}.\]
Hence $\chi_n<\mu_n,\nu_n>\chi_n^{-1}$ satisfies:
\begin{itemize}
\item
\[|\tilde{\lambda}_n|^{-1}<|z_{\chi_n\nu_n\chi_n^{-1},l}|\le|z_{\chi_n\nu_n\chi^{-1}_n,u}|<|\tilde{\lambda}_n|\]
\item
\[\lim_n\left\{\frac{(|z_{\chi_n\nu_n\chi_n^{-1},u}|+1)(|\tilde{\lambda}_n|+1)}{|\tr(\nu_n)|(|\tilde{\lambda}_n|-|z_{\chi_n\nu_n\chi_n^{-1},u}|)},
\frac{(|z_{\chi_n\nu_n\chi_n^{-1},l}|+1)(|\tilde{\lambda}_n|^{-1}+1)}{|\tr(\nu_n)|(|z_{\chi_n\nu_n\chi_n^{-1},l}|-|\tilde{\lambda}_n|^{-1})}\right\}=0\]
\end{itemize}
 conditions of Remark \ref{classical-rem-2}.
\end{proof}
Hence we have completed proof Theorem \ref{2-fixed-point}.
\end{proof}
\section{Proof of Main Theorem}
\begin{thm}\label{main-1}
There exists $\epsilon>0$ such that every $2$-generated Schottky group $\G$
with $D_\G<\epsilon$ is a classical Schottky group.
\end{thm}
\begin{proof}
This follows from Theorem \ref{t-space} and Theorem \ref{2-fixed-point}.
\end{proof}
\begin{proof}[Proof of Theorem \ref{main}]
Let $\G'$ be a non-elementary finitely generated Kleinian
group. Selberg lemma implies $\G'$ contains a torsion-free subgroup $\G''$ of
finite index, in particular $\mathfrak{D}_{\G'}=\mathfrak{D}_{\G''}$.\par
Note that if $\G''$ is geometrically infinite with $\Omega_{\G''}\not=\emptyset$ then
$D_{\G''}=2$, this implies $\mathfrak{D}_{\G''}=2$ for geometrically
infinite groups. So we can assume $\G''$ is geometrically finite when $\mathfrak{D}_{\G''}<2$.\par
If $\G''$ contains parabolic of rank $l_{\G''}$ then $D_{\G''}\ge l_{\G''}/2$.
Hence, for sufficiently small Hausdroff dimension $D_{\G''}$, we can assume $\G''$ is convex-cocompact
of second kind.\par 
It follows from Ahlfors finiteness theorem, that
$\Omega_{\G''}/{\G''}$ consists of finite number of compact Riemann surfaces.
Let $S$ be a component of $\Omega_{\G''}/{\G''}$. If $S$ is incompressible then
$\pi_1(S)$ is a surface subgroup of $\G''$. Since $1=D_{\pi_1(S)}\le D_{\G''}$, if $D_{\G''}$<1, we
may assume $S$ is compressible. So we can decompose $\G''$ along the compression
disk. After repeating the decomposition process finitely many times we are left with 
topological balls, i.e. $\mathbb{H}^3/\G''$ is a handle body. This implies $\G''$ is a finitely generated 
free purely loxodromic Kleinian group of second kind, i.e. $\G''$ is a Schottky group.\par 
By assuming the limit set have sufficiently small Hausdorff dimension 
we have reduced the general case to the case of 
Schottky groups. Now it follows from Marden's rigidity theorem, all
Schottky groups of the same rank are quasiconformally
equivalent. Therefore we have from Theorem \ref{main-1},  
there exists $\lambda>0$ such that all non-elementary finitely generated Kleinian $\G'$ 
with $\mathfrak{D}_{\G'}\le\lambda$ contains a classical Schottky group of finite index.
It follows that we have a strict lower bound on the Hausdorff dimension of all non-classical
Schottky group.
\end{proof}
\text{}\\
E-mail: yonghou@math.uic.edu
\pdfbookmark[1]{Reference}{Reference}
\bibliographystyle{plain}

\end{document}